\documentclass[a4paper,leqno]{CAIM_Sciendo_open_article_BY4}

\newcommand{\issue}{1}
\newcommand{\lastprvs}{0}
\newcommand{\DOI}{\thecwyear-XXXX}

\newcommand{\Tsteel}{T_{\mathrm{steel}}}
\newcommand{\Tflux}{T_{\mathrm{flux}}}

\newcommand{\Ksteel}{K_{\mathrm{steel}}}
\newcommand{\Ktissue}{K_{\mathrm{fabric}}}
\newcommand{\Ytissue}{\kappa_{\mathrm{fabric}}}
\newcommand{\htissue}{h_{\mathrm{fabric}}}
\newcommand{\htissuemin}{\htissue^{\min}}
\newcommand{\htissuemax}{\htissue^{\max}}
\newcommand{\wtissue}{w_{\mathrm{fabric}}}
\newcommand{\TmaxQ}{{T_{\max}^Q}}
\newcommand{\TmaxL}{{T_{\max}^L}}
\newcommand{\Cptissue}{C_{p,\textrm{fabric}}}
\newcommand{\celsius}{\textrm{\textcelsius}}

\usepackage{mathrsfs}
\usepackage[ansinew]{inputenc}
\usepackage{float}
\usepackage{afterpage}
\let\captionbox\relax
\usepackage{caption}
\DeclareCaptionLabelSeparator{dotquad}{.\quad} 
\usepackage{epstopdf}
\usepackage{blindtext}
\usepackage{fancyvrb}

\makeatletter
\renewcommand\@biblabel[1]{#1.}
\makeatother

\newcommand{\unit}[1]{{\color{blue}\mathrm{\left[#1\right]}}}

\allowdisplaybreaks

\renewcommand{\phi}{\varphi}

\newcommand{\authone}{F. Bagagiolo}
\newcommand{\authtwo}{E. Bertolazzi}
\newcommand{\auththree}{L. Marzufero}
\newcommand{\authfour}{A. Pegoretti}
\newcommand{\authfive}{D. Rigotti}
% ...

\newcommand{\authors}{\authone,~\authtwo,~\auththree, ~\authfour, ~\authfive}

\title{Modelling of the bonding process for a non-woven fabric:\\ analysis and numerics}

\author{\authone\aff{1},
        \authtwo\aff{2},
        \auththree\aff{3},
        \authfour\aff{2},
        \authfive\aff{2}}

\affiliation{\aff{1} Department of Mathematics, University of Trento, Italy
             \aff{2} Department of Industrial Engineering, University of Trento, Italy
             \aff{3} Faculty of Economics and Management, Free University of Bozen-Bolzano, Italy;\\ member of INdAM-GNAMPA research group}

\begin{document}

\setcounter{page}{\lastprvs+1}

\def\vet#1{{{\boldsymbol{#1}}}}

\maketitle

\begin{abstract}
\phantomsection\addcontentsline{toc}{section}{\numberline{}Abstract}

This paper presents research conducted at the University of Trento addressing an industrial challenge from Fater S.p.A. regarding the thermal bonding of non-woven fabrics for diaper production. The problem consists in a possible analysis of the behavior of the bonding process of a non-woven fabric. In particular, the bonding process is not given by the use of some kind of glue, but just by the pressure of two fiber webs through two high-velocity steel-made rollers. The research comprised the formulation and theoretical as well as numerical analysis of analytical, mechanical and thermal models for the stress-strain behavior of the non-woven fabric's fibers and for the bonding process with heating effects.
\end{abstract}


\AMScode{65L05; 65M06; 80A19; 80M20.}

\section{Introduction}
A non-woven fabric is a manufactured material made from polymer fibers which are bonded together through industrial processes that usually involve pressure and heating effects. There are several purposes for the use of nonwoven fabrics, such as medical fabrics, filters, or geotextile ones. Among them, we can also find the use for the production of diaper linings. The paper is indeed concerned with the modelization of a bonding process for a non-woven fabric used for the manufacturing of diapers. In particular, the problem was posed by the Italian company Fater S.p.A, specialized in the market of absorbent products, with the aim of a better understanding of their industrial bonding process for producing diapers, in the event that the nonwoven fabric is provided with different physical and mechanical parameters. Such a bonding process performed by the company does not require the use of some kind of glue but involves just the pressure of two fabric webs passing through two high-velocity steel-made rollers that compress them. In this way, an instantaneous heat diffusion and local fusing effect take place and the two webs of non-woven fabric will get glued. The mechanics behind this kind of bonding process is also partly shared by the bonding process used to produce a non-woven fabric; see \cite{Patent1, Patent2} for more details. Our purpose is then to model the bonding process, both from an analytical and numerical point of view, as a function of the compression and the velocity of the rollers. To this aim, two different models for the heat diffusion were studied, a linear and a parabolic one. The first is less accurate and involves just ordinary linear differential equations for temperature (see \cite{Cengel}), whereas the second aims to enhance the compression of the fusing process by exploiting the heat equation (see \cite{Formaggia}).
\par
As regards the experimental part, the use of laboratories and equipments, we made use of the expertise of some of the authors who have already worked on similar models as in \cite{Pegoretti}. See also the references therein for the industrial importance of such materials.
\par
This allows us to study the behavior of the resulting diaper lining even in the case some parameters like thickness are different. Our investigation is detailed as in the following sections.
\par
In Section \ref{description_nonwoven_fabric}, we briefly describe the non-woven fabric, the way it is usually produced and the intended use. Some pictures (made by an electronic microscope) of the non-woven fabric used by the company are reported. In particular, such pictures clearly enlighten the fact that the bonding process exploited for the production of the fabric presents an evident thermal effect. The following chapters will then study from an analytical, mechanical and numerical point of view a new model taking into account those thermal effects.
\par
In Section \ref{polypropylene}, after our laboratory experiments, we get that the thermal effect, i.e. the temperature increment of the non-woven fabric, is due to the almost instantaneous compression process. We then collect other experimental results as well as data from literature in order to precisely fit the nature of the fabric and estimate the involved parameters, such as specific heat, Young's modulus, etc. Thereafter, some analytical models inspired by the so-called heat equation are formulated and studied.
\par
In Section \ref{bondingprocess}, we further specify the analytical model in Section \ref{polypropylene}, making it more suitable for the real process of compression and bonding given by the passage of the fabric webs through the rollers. Then a one-dimensional and a bi-dimensional numerical model and the relative numerical results are reported. In particular, such results clearly show the patterns of effective thermally bonding as function of the compression of the fabric and of the velocity of the assembly line.
\par

In Section~\ref{conclusions}, we outline possible future developments and make some further comments on the interpretation and applicability of the results.
\section{Description of the non-woven fabric}
\label{description_nonwoven_fabric}
A non-woven fabric is a fabric-like material produced by bonding together staple short or continuous long fibers through chemical or mechanical processes. Many treatments for producing non-woven fabrics exist. For example, one can apply heat and pressure for bonding at limited areas of a non-woven film by passing it through the nip between heated calender rolls either or both of which may have patterns of lands and depressions on their surfaces. During such a bonding process, depending of the types of fibers making up the non-woven film, the bonded regions may be formed independently of external influence or aid, i.e. the fibers of the film are melt fused at least in the pattern areas or with the addition of an adhesive. The advantages of thermally bonded non-woven fabrics include low energy costs and speed of production. 
\par
Non-woven fabrics can also be made by other processes. For more details see the patents \cite{Patent1} and \cite {Patent2} and references (and other patents) therein. However, in all of the non-woven fabrics, the producing process usually realizes punched rigid spots on the pattern giving them a sort of frame. For simplicity, from now on such punched rigid spots will be called \textit{tags}.
\par
The non-woven fabric provided by the company is used for producing diapers linings and is formed by polypropylene. Generally, non-woven fabrics are also used for medical and sanitary stuff like hospital gowns, wipes and masks and may be also formed by other types of polymers, like polyolefin and polyester. In particular, the company exploits a further bonding process to thermally bond pieces of raw non-woven fabric (made in turn by a bonding process like the ones described above) in order to build pieces of a diaper lining. The technical equipment used for such a bonding process is patented by the company and very similar to the one in \cite{Patent2}: two fiber webs pass through a pair of high velocity steel-made rollers which thermally fuse them in the pattern areas by applying a suitable almost instantaneous pressure. In principle, no heating process is needed for the rollers (for more details, see the Section \ref{bondingprocess}). Therefore, the description of the production suggests the numerical modelling of the bonding process.

\subsection{First collection of pictures}
\begin{figure}[H]
  \begin{center}
    \includegraphics[scale=0.20,angle=0]{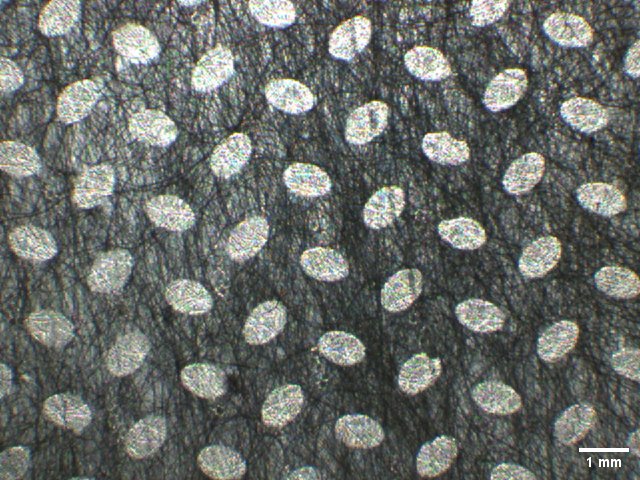}
    \includegraphics[scale=0.20,angle=0]{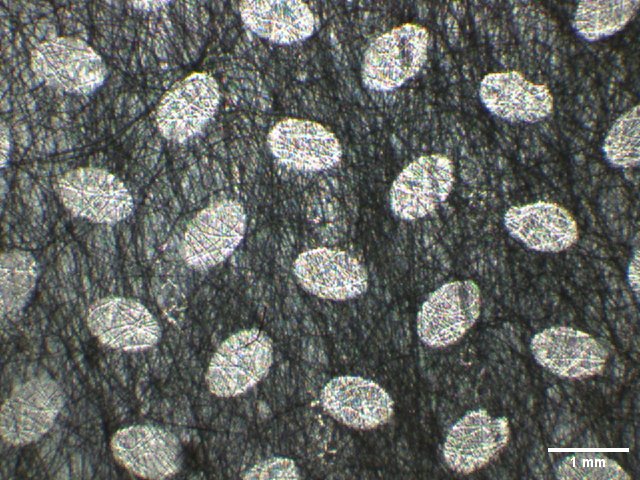}
    \includegraphics[scale=0.20,angle=0]{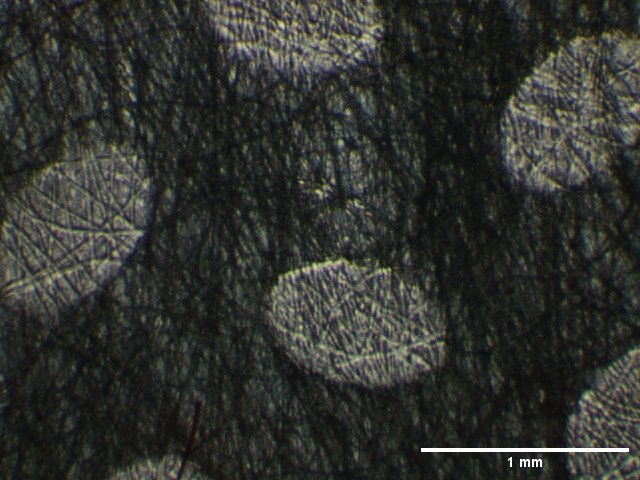}
  \end{center}
  \caption{First collection of pictures illustrating the elliptic shape of the tags and randomness of fibers.}
  \label{fig:microscopio1}
\end{figure}

In this section, we report some pictures of pieces of nonwoven fabric provided by the company, made by the Heerbrugg Wild M3Z optical microscope (see Figure~\ref{fig:microscopio1}), without any focus on the thermally bonded pattern. The presence of tags and fibers is shown anyway. 

\subsection{Second collection of pictures}
\begin{figure}[H]
  \begin{center}
    \includegraphics[scale=0.15]{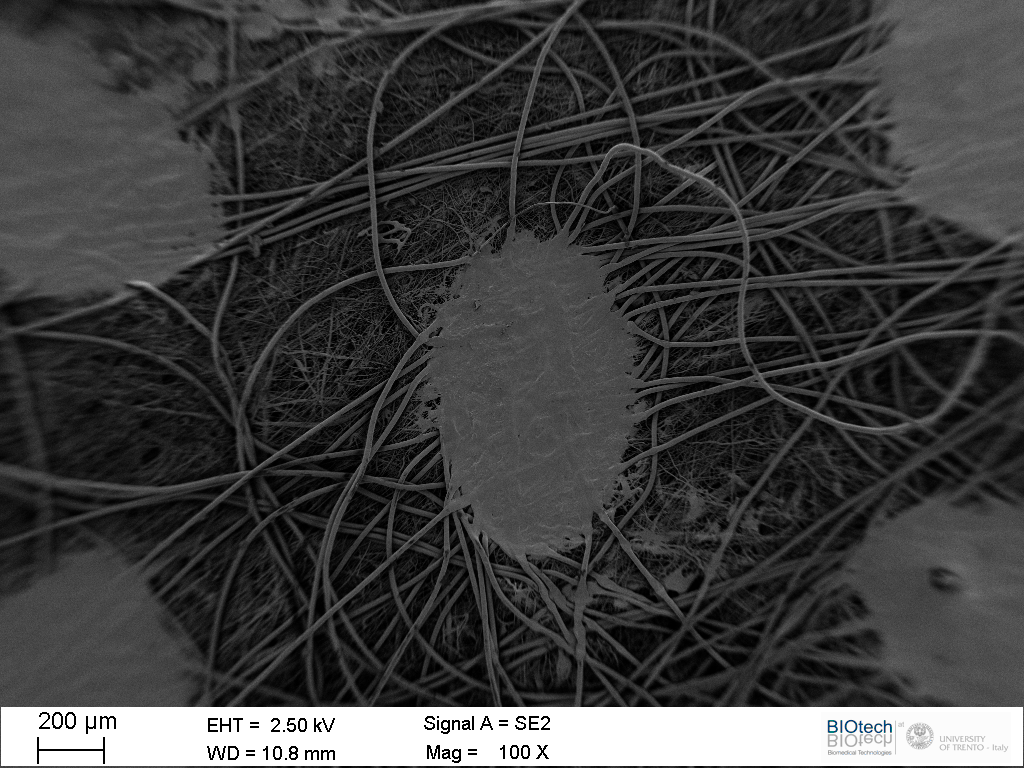}
    \includegraphics[scale=0.15]{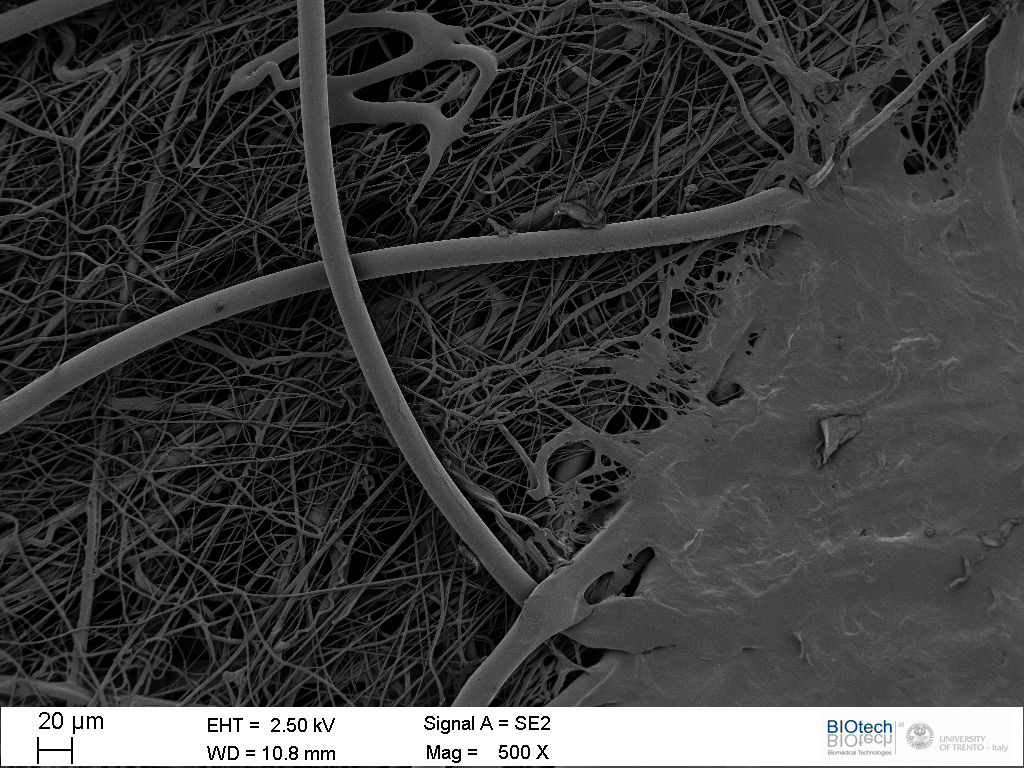}
    \includegraphics[scale=0.15]{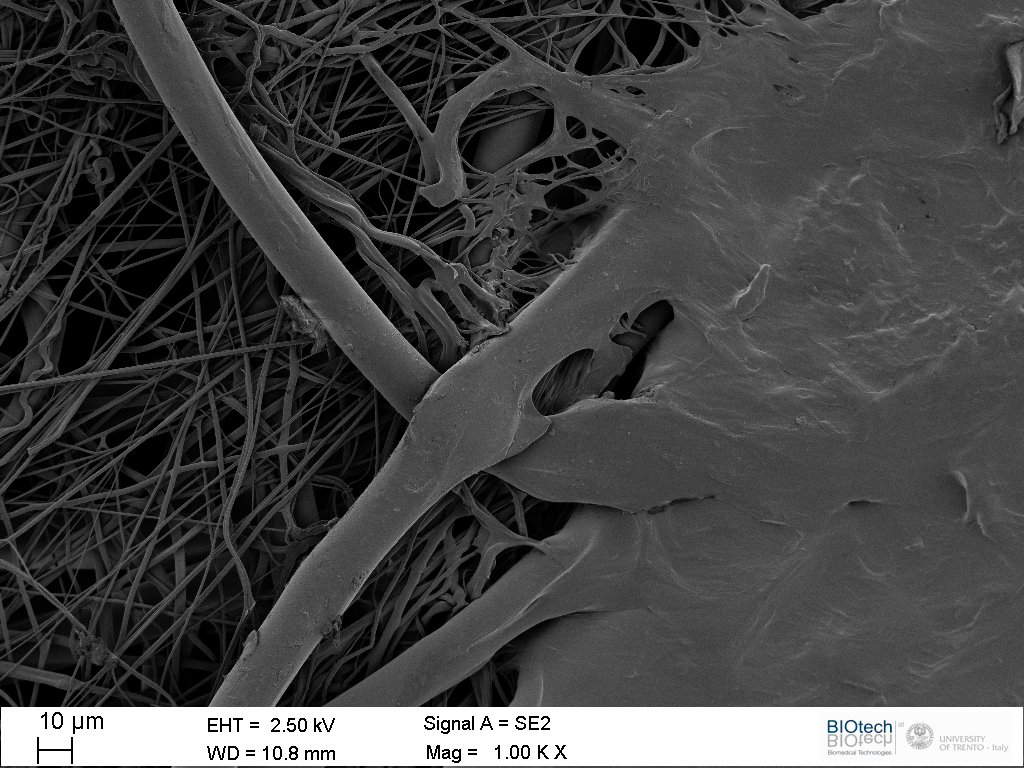}
  \end{center}
  \caption{Second collection of pictures illustrating the fused pattern areas near elliptic tags induced by the thermally bonding process.}
  \label{fig:microscopio3}
\end{figure}

Here we report another collection of pictures (see Figure~\ref{fig:microscopio3}) with a detailed focus on the thermally bonded areas, that unequivocally shows that, around the elliptic tags, the almost instantaneous compression process induces a thermal bonding with a resulting local fusing. The microscope we used is Zeiss Supra 40 field emission scanning electron microscope (FESEM), operating at an accelerating voltage of $2.5\ \unit{kV}$. A thin platinum palladium conductive coating was deposited on the surface of the samples before the observations.

\section{Polypropylene data and modelling}
\label{polypropylene}
The aim of this section is to compute the temperature increment and the temperature decay of the non-woven fabric through a suitable compression rate. We will provide two models, a linear and a quadratic one. To do that, we will need some data of the fabric, polypropylene and steel, which we will recover a little from the literature and a little from laboratories experiments.
\subsection{Tacticity and some data from literature}
In order to establish which is the polypropylene's type (atactic, isotactic, syndiotactic) of the non-woven fabric, we performed a spectrum analysis, made with a Fourier transformed infrared (FT-IR) spectrophotometer, using a Avatar 330 by the Thermo Fisher company. The spectrum, as in figure~\ref{fig:spettro}, is the result of $64$ scans with a resolution of $4\ \unit{cm}^{-1}$, and turned out that the polypropylene is atactic.


\begin{figure}[!htb]
  \begin{center}
    \includegraphics[scale=0.35]{figure/polipropilene}
  \end{center}
  \caption{Comparative analysis of the non-woven fabric and reference atactic polypropylene samples. The spectra show matching characteristic absorption peaks}
  \label{fig:spettro}
\end{figure}

\begin{figure}[!htb]
  \begin{center}
    \includegraphics[scale=0.4]{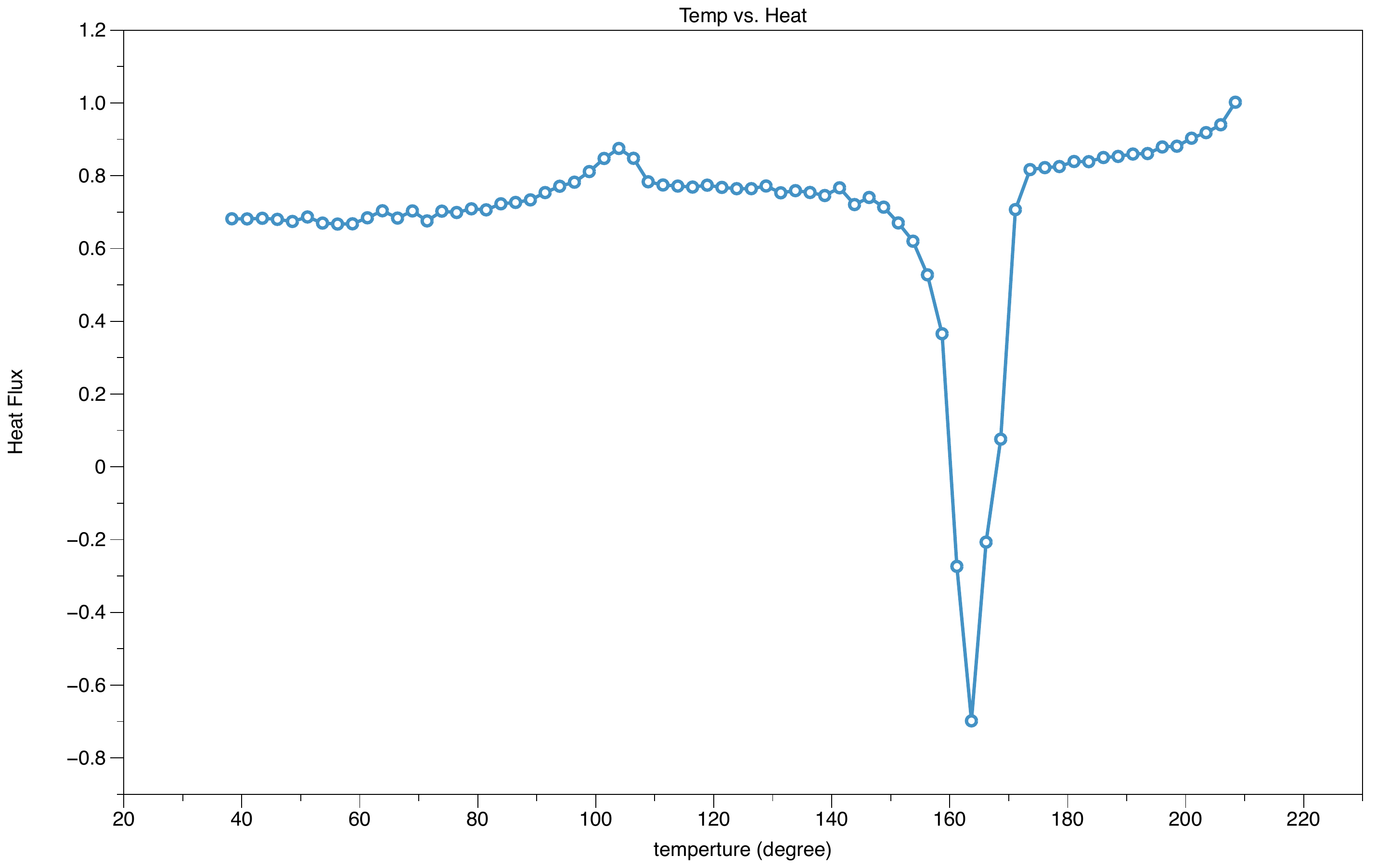}
  \end{center}
  \caption{Differential Scanning Calorimetry showing the endothermic melting transition of polypropylene near $160$ $\unit{\celsius}$. The negative heat flow peak corresponds to the energy absorption during fusion phase}
  \label{fig:fusion}
\end{figure}

For the description and the numerical modelling of the bonding process, we need some polypropylene data, that we collect in the following Table~\ref{tab:polipropilene} (see \cite{Passaglia, Maier, Young,Gianotti1968}):
\begin{table}[!htb]
\centering
  \caption{Physical parameters of polypropylene}
   \label{tab:polipropilene}
  \begin{tabular}{llll}
    Parameter     & Value & Unit \\
    \hline
    Melting point        & $130$--$171$   & $\unit{\celsius}$        & measured $160$ $\celsius$ \\
    Density              & $855$--$946$   & $\unit{kg/m^3}$   & atactic: $866$ \\
    Thermal conductivity & $0.17$--$0.22$  & $\unit{W/(m K)}$  &  at $23$ $\celsius$ typ $0.17$ \\
    Specific heat        & $581.97$--$2884.71$ & $\unit{J/(kg\, K)}$ & temp. in K \\
    Young's module         & $1300$--$1800$ & $\unit{N/(mm)^2}=\unit{MPa}$               & \\
    \hline
  \end{tabular}
\end{table}
Due to the high variability of the fusion temperature $130$--$171 \ \unit{\celsius}$, an accurate measurement is performed fixing at $160\ \unit{\celsius}$ (see figure~\ref{fig:fusion}) the fusion temperature. In order to detect that, polypropylene is heated slowly with a constant temperature rate. The heat flux to the polypropylene is measured and, when a negative peak of heat is detected, it means that polypropylene is changing phase from solid to liquid.

Since the rollers which thermally bond the pieces of non-woven fabric are steel-made, similarly to the polypropylene, we need some physical parameters of the steel, that we collect in the following Table~\ref{tab:acciaio}~(see \cite{steel3, steel4, steel2, steel1}):

\begin{table}[!htb]
  \caption{Physical parameters of steel}
  \label{tab:acciaio}
  \begin{center}
  \begin{tabular}{lll}
    Parameter     & Value & Unit \\
    \hline
    Melting point               & $1400$--$1530$ & $\unit{\celsius}$                   \\
    Density                     & $7500$--$8000$ & $\unit{kg/m^3}$   \\
    Thermal conductivity (stainless) & $15$--$18$     & $\unit{W/(m K)}$  \\
    Thermal conductivity        & $44$--$80$     & $\unit{W/(m K)}$  \\
    Specific heat               & $500$          & $\unit{J/(kg K)}$ \\
    Young's modulus              & $180000$       & $\unit{N/(mm)^2}=\unit{MPa}$ \\
    \hline
  \end{tabular}
  \end{center}
\end{table}




\subsection{Determination of fabric weight by unit area}
In this section, we report the density of our non-woven fabric, determined by a $15\ \unit{cm} \times 15\ \unit{cm}$ fabric sample
(see figure~\ref{fig:weight}). The next figure shows the fabric piece and its weight, determined by a precision scale Netzsch DSC204.
\begin{figure}[!htb]
  \begin{center}
    \includegraphics[height=4.7cm]{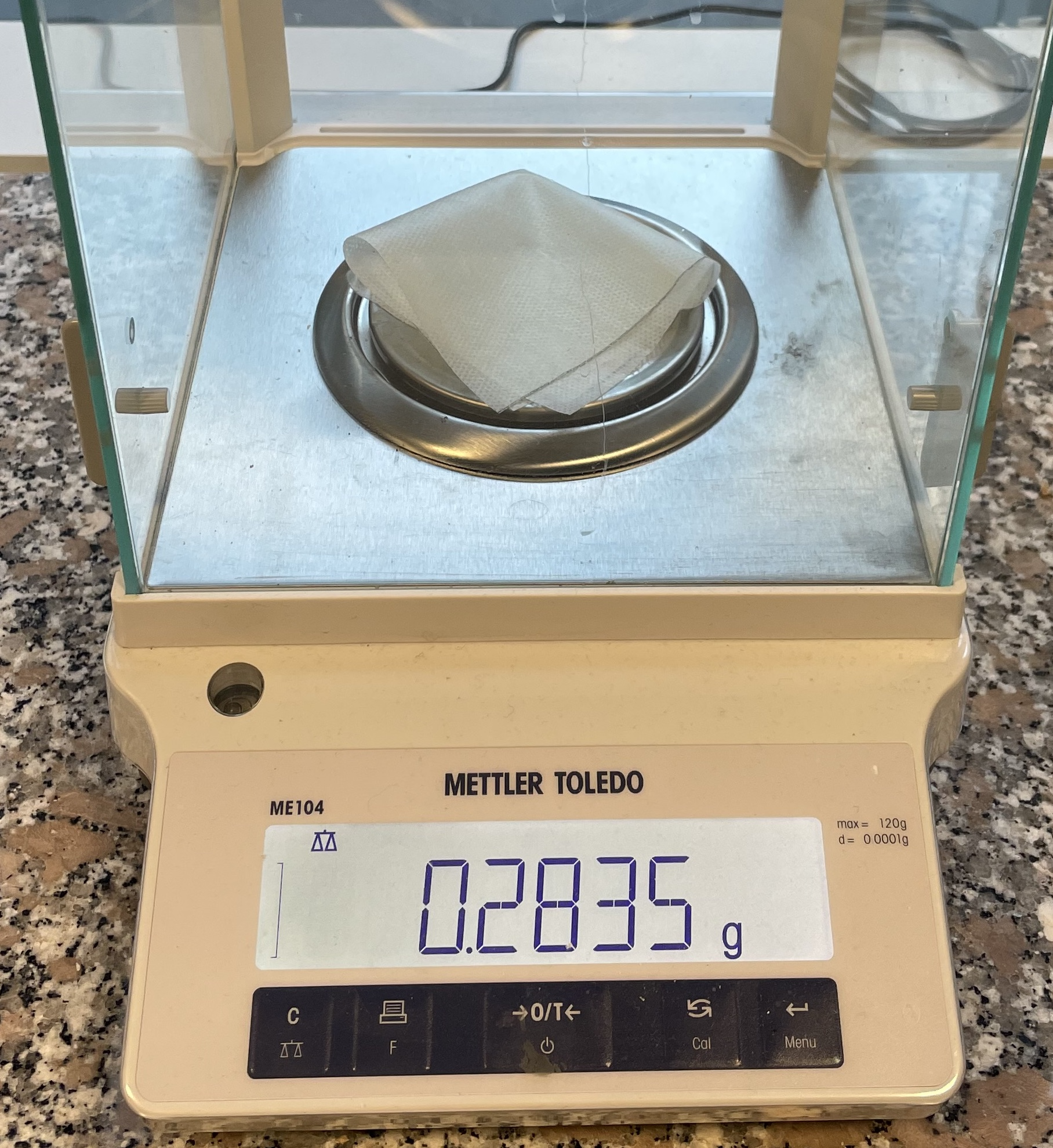}
    \includegraphics[height=4.7cm]{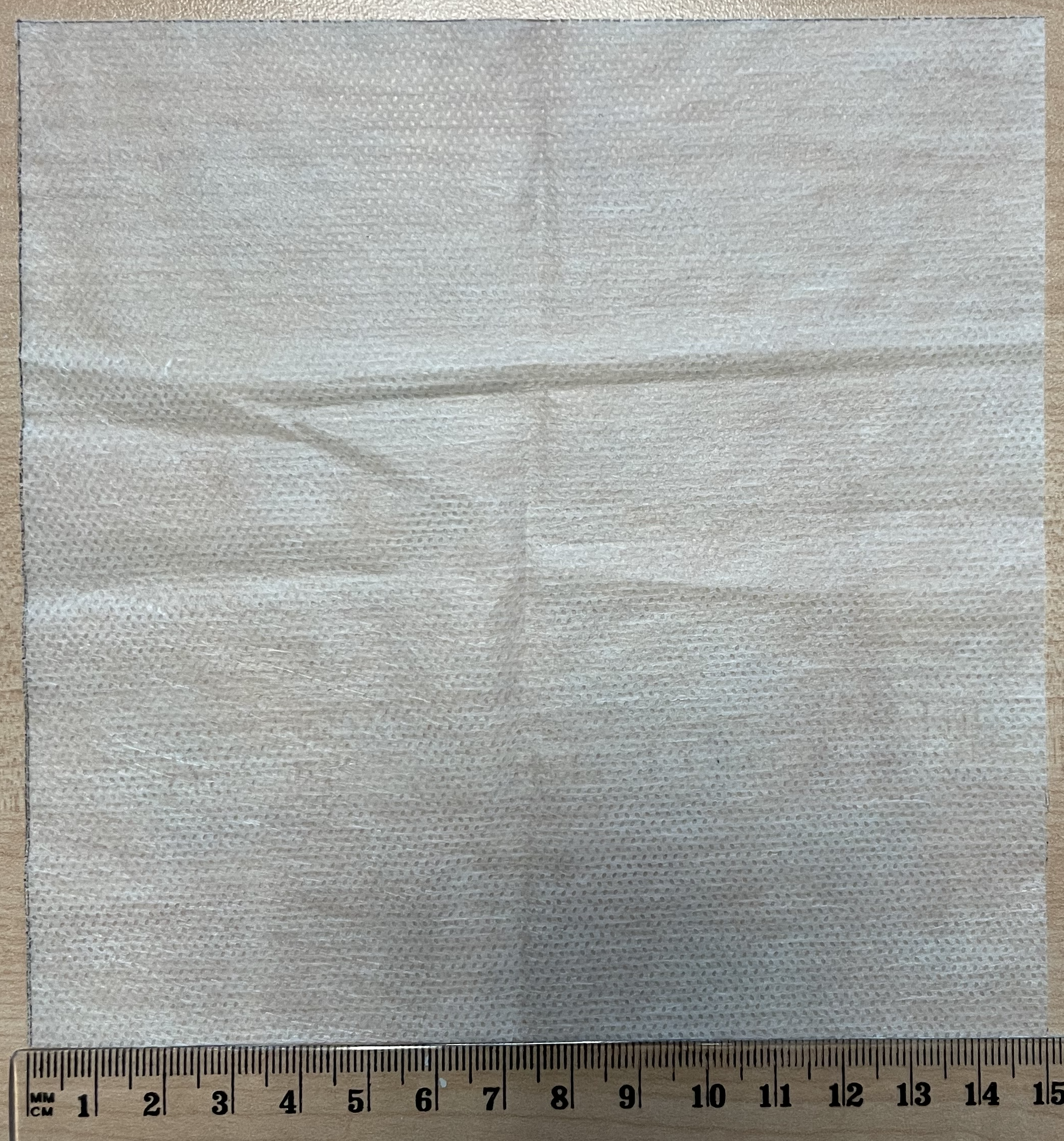}
    \includegraphics[height=4.7cm]{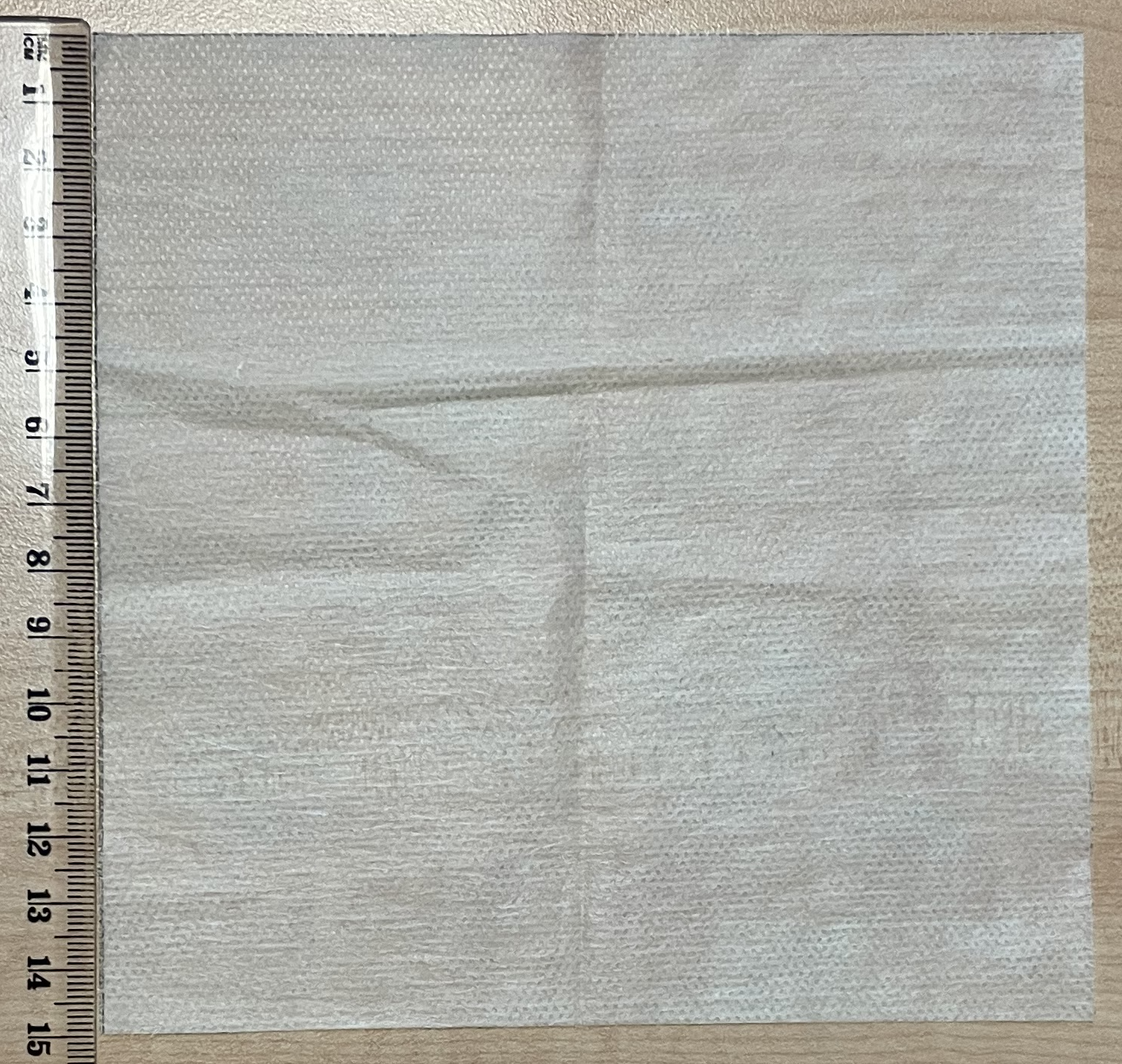}
  \end{center}
  \caption{On the left the precision scale used for the weight determination and on the right the non-woven fabric piece used for the experiment.}
  \label{fig:weight}
\end{figure}
It follows that the weight by unit area $\wtissue$ is:
\begin{equation}\label{eq:wT}
   \wtissue=\dfrac{[\textrm{mass}]}{[\textrm{area}]}=
   \dfrac{0.2835 \ \unit{g}}
         {150 \ \unit{mm}\cdot 150 \ \unit{mm}} = 0.0126 \ 
   \unit{\dfrac{kg}{m^2}}.
\end{equation}
As long as the nonwoven fabric is fully compressed, from the density of the polypropylene (see Table \ref{tab:polipropilene}) we can infer its thickness:
\[
   \htissuemin = \dfrac{[\textrm{weight}]/[\textrm{area}]}{[\textrm{density}]}
     = \dfrac{0.0126 \ \unit{\dfrac{kg}{m^2}}}
             {900 \ \unit{\dfrac{kg}{m^3}}}
     = 0.000014\,\unit{m}
     = 0.014\,\unit{mm}
     = 14\,\unit{\mu m}.
\]
\subsection{Computation of the thickness: pressure and displacement}

In this subsection, we want to establish the thickness of the 
not compressed non-woven fabric. To do that, we use an indirect approach based on experimental measurements made with the dynamometer Instron 5969, at a speed of $1\ \unit{mm}/\unit{min}$, of the laboratories in the Department of Materials Engineering at the University of Trento.

\begin{figure}[!htb]
  \begin{center}
    \includegraphics[width=12cm]{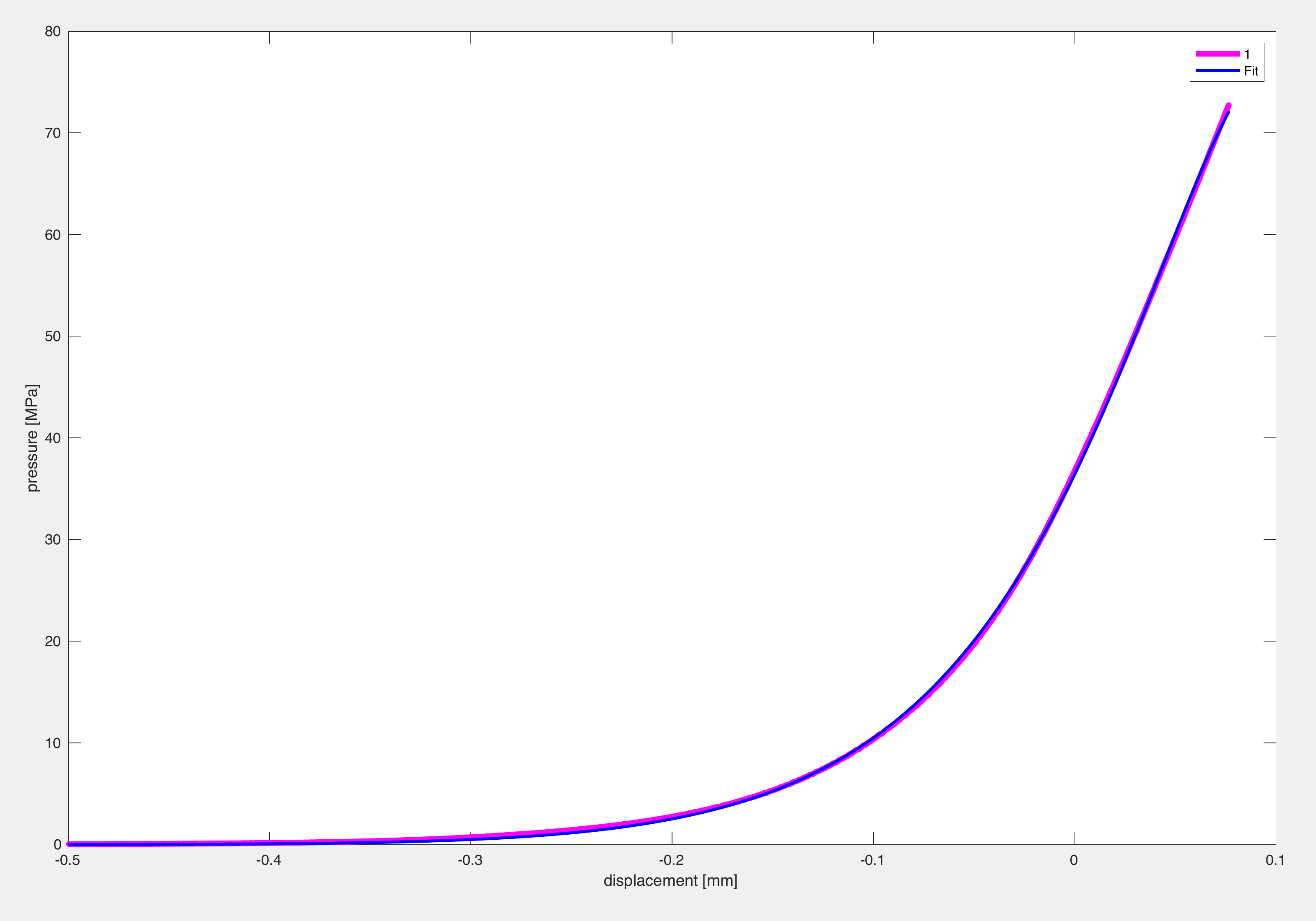}
  \end{center}
  \caption{Displacement/pressure graph with ten non-woven fabric sheets}
  \label{fig:10:sheets}
\end{figure}

In particular, at first we characterize the displacement/pressure curve of the dynamometer with the press disk made by seven circle tags without the non-woven fabric. In this way, we determine the offset corresponding to zero thickness to be eliminated in the measurement of the curve. We did two experiments and the results are the following.
\begin{figure}[!htb]
  \begin{center}
    \includegraphics[width=7cm]{figure/pressione-spostamento-free-a.pdf}
    \includegraphics[width=7cm]{figure/pressione-spostamento-free-b.pdf}
  \end{center}
  \caption{%
  Baseline displacement-pressure curve of the press system without the non-woven fabric. The plot characterizes the intrinsic mechanical compliance of the press apparatus, which will be subtracted from subsequent measurements with fabric samples to isolate the material's deformation behavior}
  \label{fig:press:free}
\end{figure}
The fitting of the pressure is the following (see figure~\ref{fig:press:free})
\[
  P_{\textrm{base}}(x) =
  \dfrac{5461.352911\,\max(0,x)^{2.92599}}
  {0.0038158166 + 6.4490865\,\max(0,x)^{1.624481}}.
\]
After this measurement, ten sheets of non-woven fabric are put under the press disk, getting the following pressure curve
of figure~\ref{fig:10:sheets}.
Fitting again the data, we obtain
\[
  P_{\textrm{base}+\textrm{10 fabric}}(x) =
  \dfrac{716.33893 \max(0,x+0.9703)^{12.67189680}}
        {14.10752 + 0.92219399\max(0,x+0.9703)^{30.037944}}.
\]
This shows that the pressure starts growing at $x=-0.9703$ which can be assumed as the thickness of the ten non-woven fabric sheets, and then we can set it to
\[
   \htissuemax = \dfrac{0.97\,\unit{mm}}{10} = 97\,\unit{\mu m}.
\]
Due to the very low speed of the displacement movement, the pressures are in equilibrium as follows:
\[
  P_{\textrm{base}+\textrm{10 fabric}}(x) =
  P_{\textrm{base}}(z) = P_{\textrm{base}}( x- w) =
  P_{\textrm{10 fabric}}(w),
\]
where the total displacement $x=z+w$ is the sum of the displacement of the ten non-woven fabric sheets ($w$) and the press displacement ($z$). To obtain $P_{\textrm{10 fabric}}(w)$, first of all we have to determine $x$ as a function of $w$ by solving
\begin{equation}\label{eq:xw:solve}
  P_{\textrm{base}+\textrm{10 fabric}}(x)-
  P_{\textrm{base}}(x-w)= 0.
\end{equation}
This equation can be solved numerically with respect to $x$ and assuming $x(w)$ as known. Then
\[
  P_{\textrm{10 fabric}}(w) = P_{\textrm{base}+\textrm{10 fabric}}(x(w)).
\]
The function $x(w)$ is very ill-conditioned and defined for about $w\leq-0.08$. Thus it is better to compute the inverse function $w(x)$ by solving~\eqref{eq:xw:solve} w.r.t. $w$. Plotting $w(x)/10$ and $P_{\textrm{base}+\textrm{10 fabric}}(x)$, we obtain
\begin{figure}[!htb]
  \begin{center}
    \begin{tabular}{cc}
      \vtop{\null\hbox{\includegraphics[height=8cm]{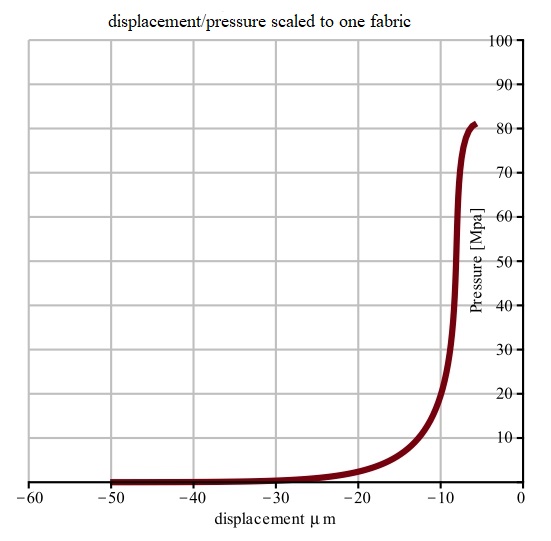}}} &
      \vtop{\null\hbox{\includegraphics[height=9cm]{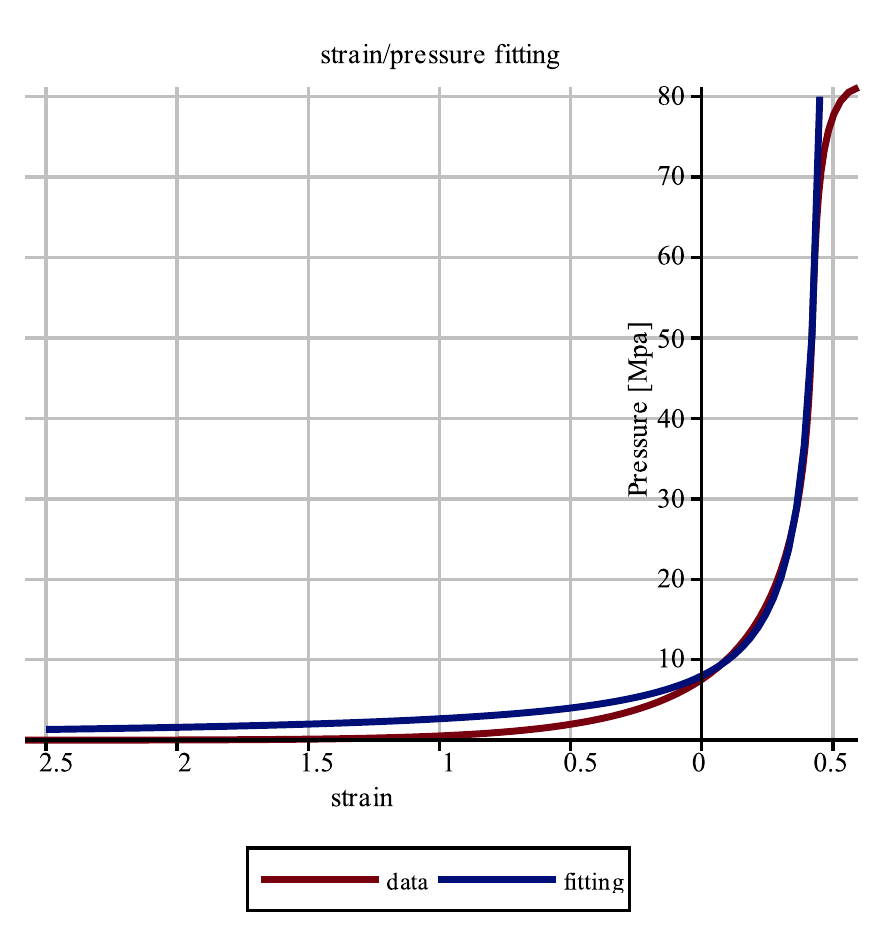}}}
    \end{tabular}
  \end{center}
  \caption{On the left: displacement/pressure graph of a single sheet of non-woven fabric\\
  On the right: strain/pressure curve fitting}
  \label{fig:pp:pressure:fitting}
\end{figure}
and this is an approximation of the displacement/pressure curve on a single sheet of non-woven fabric. When the displacement is $0$, the fabric has non physical thickness $0$ and graph stops at about $-6 \ \unit{\mu m}$. The graph pressure is negligible until displacement is $-30 \ \unit{\mu m}$, and thus we can correct the thickness of the fabric (under mild compression) to $30 \ \unit{\mu m}$.
\par
Hence the fabric under compression at about $5 \ \unit{MPa}$ is fully compressed (the fibers are compacted) and the thickness is about $14\,\unit{\mu m}$. When pressure grows until $80\,\unit{MPa}$, the thickness is reduced to about $7\,\unit{\mu m}$, and the non-woven fabric assumed as a single block is under a strain of $7/14=0.5$. This information will be used to compute the increment of the temperature of the fabric under fast compression.
\par
To make computations workable, we derive from Figure~\ref{fig:pp:pressure:fitting} left an approximated law linking the strain of the non-woven fabric with the pressure. In Figure~\ref{fig:pp:pressure:fitting} right, a plot of the strain calculated with respect to the compressed fabric and the pressure is fitted with a hyperbolic curve. The $x$-axis contains the strain:
\[
   [\textrm{strain}]
   = \dfrac{[\textrm{displacement}]+\htissuemin}{\htissuemin}
   = \dfrac{[\textrm{displacement}]+14 \ 
   \unit{\mu m}}{14 \ \unit{\mu m}}.
\]
Zero strain means that the fabric is fully compressed and positive strain means that the fabric is compressing. Negative strain means that the fabric is not fully compressed and the pressure is relatively low.

As we can notice, the pressure in the negative strain part is overestimated, while the positive part is well captured up to a strain of $0.4$ (compression of $40\%$). However, we are interested in the estimation of the heating of the non-woven fabric, and this part of the curve produces a very low effect. The approximated/fitted Young's modulus becomes
\begin{equation}\label{eq:Y:module:tissue}
  \Ytissue(s)
  = \dfrac{16}{1-2s}\ \unit{MPa}
  = \dfrac{16\cdot 10^6}{1-2s}\ \unit{Pa}.
\end{equation}

\subsection{Computation of the temperature increment}
\begin{figure}[!htb]
  \begin{center}
    \includegraphics[scale=0.5]{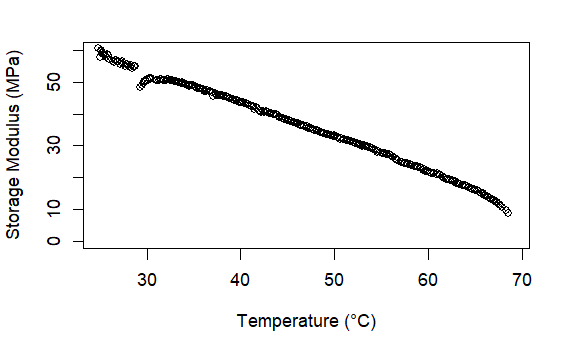}
  \end{center}
  \caption{Storage modulus/temperatue plot of the non-woven fabric}
  \label{fig:stiffness:temperature}
\end{figure}

From the following experimental figure~\ref{fig:stiffness:temperature} obtained at the University of Trento, it follows that $\Ytissue$ (i.e., the Young or Storage modulus) of the non-woven fabric linearly decreases by temperature. Dynamic mechanical analysis (DMA) tests were carried out using a TA Instrument DMA Q800 device, in the temperature range from $25\ \unit{\celsius}$ to $120\ \unit{\celsius}$, with a heating rate of $3\ \unit{\celsius/(min)}$, a strain amplitude of $0.05\%$ and a frequency of $1\ \unit{Hz}$. Through this analysis, it was possible to evaluate the Storage modulus, the loss modulus and the loss tangent as a function of the temperature. 
In particular, the Young modulus linearly decreases from $T_0=20$ $\unit{\celsius}$ and reaches $0$ at about $90$ $\unit{\celsius}$. The value of the modulus of the non-woven fabric is not constant and depends on the displacement, but we can assume that its value linearly decreases by temperature in any condition. Thus we set
\[
  \kappa(s,T)
  = \Ytissue(s) \max\left(0,\dfrac{\TmaxL-T}{\TmaxL-T_0}\right).
\]
Following \cite{Young}, instead of a linear decreasing, a parabolic fitting with $\kappa(T)=0$ near the fusion temperature ($\approx 160$ $\unit{\celsius}$) can be used. Moreover, it suggests an exponential fitting for Storage modulus depending on temperature. However we used linear and parabolic fitting to maintain the model as simple as possible.
Therefore
$$
   \kappa(s,T)
   = \Ytissue(s)\max\left(0,\dfrac{\TmaxQ-T}{\TmaxQ-T_0}\right)^2.
$$
Using $\kappa(s,T)$, the force and the pressure per unit area become
\begin{equation}\label{eq:pressure:tissue}
  P(s,T) = \dfrac{F(s,T)}{[\textrm{area}]} = \kappa(s,T)s.
\end{equation}
The infinitesimal work done on the non-woven fabric (at fixed temperature) is
\[
   \mathrm{d}W(s,T) = F(s,T) \mathrm{d}s,
\]
and the variation of temperature is done by the work on the fabric as follows:
\[
    \dfrac{\mathrm{d} T(s)}{\mathrm{d}s} =
    \dfrac{\Delta[\textrm{work}]'}{\Cptissue[\textrm{mass}]}
    =
    \dfrac{\dfrac{\Delta[\textrm{work}]'}{[\textrm{area}]}}
          {C_{p,\textrm{fabric}}\dfrac{[\textrm{mass}]}{[\textrm{area}]}}
    =
    \dfrac{F(s,T)/[\textrm{area}]}{\Cptissue\cdot \wtissue}
    =
    \dfrac{P(s,T)}{\Cptissue\cdot \wtissue}
    =
    \dfrac{\kappa(s,T)s}{\Cptissue\cdot \wtissue},
\]
where $\Cptissue$ is the specific heat of the fabric. Thus, using~(\ref{eq:pressure:tissue}), assuming that the cooling temperature negligible, the temperature increment is obtained by solving the ODE
\begin{equation}\label{ode:lin}
   \dfrac{\mathrm{d} T(s)}{\mathrm{d}s} =
   \dfrac{s\,\Ytissue(s)}{\Cptissue\cdot \wtissue}\cdot \max\left(0,\dfrac{\TmaxL-T(s)}{\TmaxL-T_0}\right)
\end{equation}
for the simple linear model, and
\begin{equation}\label{ode:para}
   \dfrac{\mathrm{d} T(s)}{\mathrm{d}s}  =
   \dfrac{s\,\Ytissue(s)}{\Cptissue\cdot \wtissue}\cdot
   \max\left(0,\dfrac{\TmaxQ-T(s)}{\TmaxQ-T_0}\right)^2
\end{equation}
for the parabolic interpolation model. Using the value in the Table~\ref{tab:numeric:parameters}.
\begin{figure}[!htb]
  \begin{center}
    \includegraphics[width=8cm]{figure/heating-by-pressure.pdf}
  \end{center}
  \caption{Heating by pressure (adiabatic case)}
  \label{fig:heating1}
\end{figure}
the heating (numerically) computed is plotted in Figure~\ref{fig:heating1}
and a moderate strain of the compressed non-woven fabric is enough to heat up to the melting temperature.

\subsection{Evaluation of the temperature decay by flux}
\begin{figure}[!htb]
    \begin{center}
      \includegraphics[width=5cm]{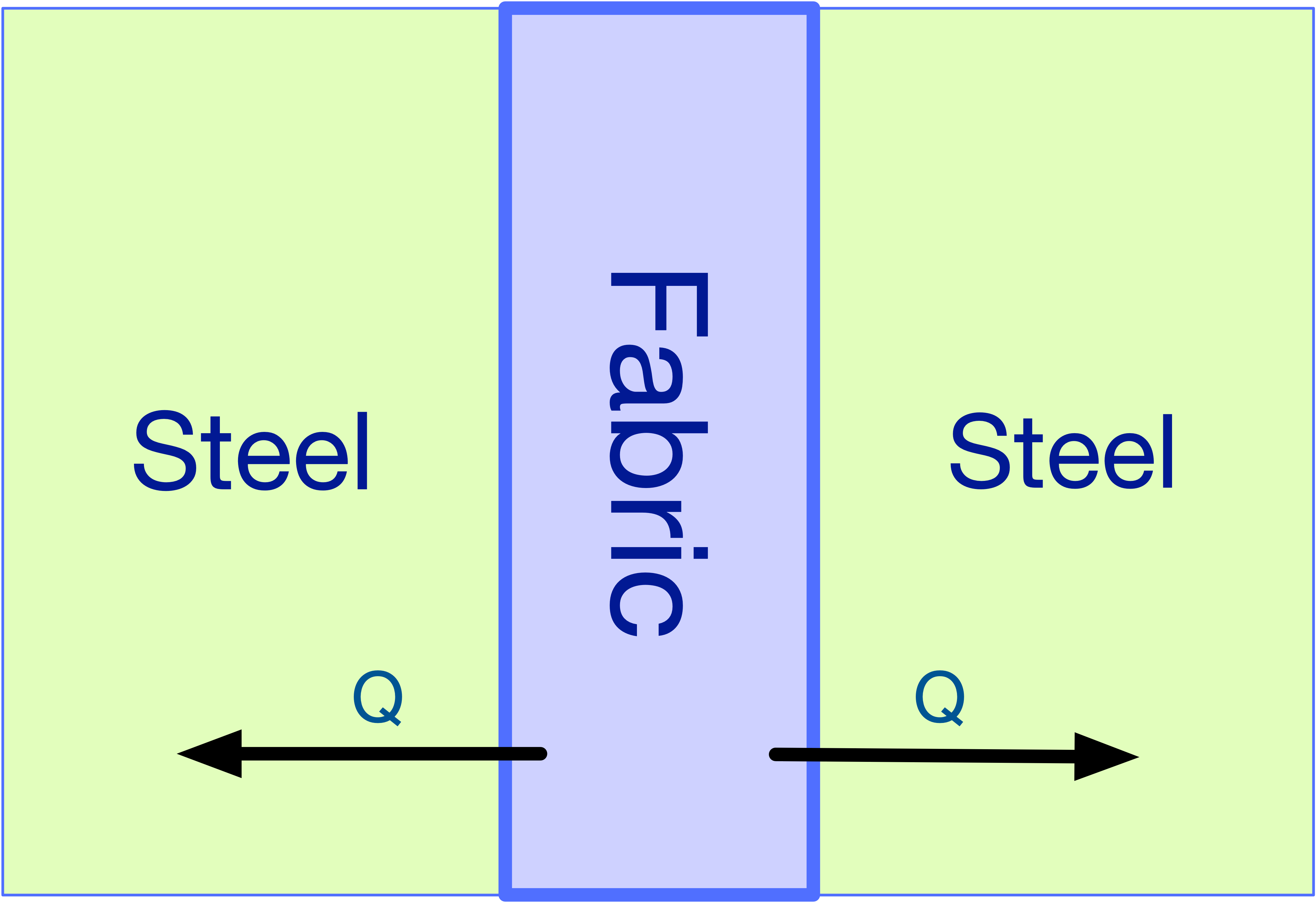}
    \end{center}
    \caption{Heat flux schematic}
    \label{fig:heat-flux}
\end{figure}
From Figure~\ref{fig:heat-flux}, assuming that the steel has fixed temperature $\Tsteel$, the thermal flux from the non-woven fabric to the steel using the Newton's Law of cooling is
\[
   \mathcal{Q} = 2\Ksteel
   \dfrac{\Tsteel-T(t)}{\htissuemin(1-s)/2} [\textrm{area}]
   \quad \unit{W},
\]
where $\Ksteel$ is the thermal conductivity of the steel. The thermal conductivity and the heat capacity of the steel are very high compared to the polypropylene ones, so that they can be assumed as infinity. The temperature variation due to the heat flux using~\eqref{eq:wT} becomes
\[
  \Cptissue[\textrm{mass}]\dfrac{\mathrm{d}\Tflux(t)}{\mathrm{d}t}
  = \mathcal{Q} = \dfrac{4\Ksteel(\Tsteel-T(t))}{\htissuemin(1-s)}  [\textrm{area}]
\]
so that
\begin{equation}\label{eq:T:cool:by:steel}
  \dfrac{\mathrm{d}\Tflux(t)}{\mathrm{d}t}
  =
  \dfrac{4\Ksteel(\Tsteel-\Tflux(t))}{\htissuemin(1-s)\Cptissue[\textrm{mass}]/[\textrm{area}]}
  =
  \dfrac{4\Ksteel(\Tsteel-\Tflux(t))}{\htissuemin(1-s)\Cptissue\wtissue}
\end{equation}
assuming that the compression is performed at a constant speed $v\,\unit{m/s}$. Then we can write
and
\begin{equation}\label{eq:Tflux}
  s = \dfrac{v\cdot t}{\htissuemin}
  \qquad\textrm{and}\qquad
  \dfrac{\mathrm{d}\Tflux(s)}{\mathrm{d}s}\dfrac{v}{\htissuemin}=
  \dfrac{4\Ksteel(\Tsteel-\Tflux(s))}{\htissuemin(1-s)\Cptissue\wtissue}.
\end{equation}
Using~\eqref{eq:Tflux} with~\eqref{ode:lin} or~\eqref{ode:para}, we have the two complete models:
\begin{equation}\label{ode:lin:full}
   \dfrac{\mathrm{d}T(s)}{\mathrm{d}s} = \dfrac{s\,\Ytissue(s)\max\left(0, \dfrac{\TmaxL-T(s)}{\TmaxL-T_0}\right)
   +\dfrac{4\Ksteel}{v(1-s)}
  (\Tsteel-T(s))}{\Cptissue\cdot \wtissue}
\end{equation}
for the simple linear model, and
\begin{equation}\label{ode:para:full}
   \dfrac{\mathrm{d}T(s)}{\mathrm{d}s} = \dfrac{s\,\Ytissue(s)\cdot
   \max\left(0,\dfrac{\TmaxQ-T(s)}{\TmaxQ-T_0}\right)^2
   +\dfrac{4\Ksteel}{v(1-s)}
  (\Tsteel-T(s))}{\Cptissue\cdot \wtissue}
\end{equation}
for the quadratic one.
We can solve numerically~\eqref{ode:lin:full} and~\eqref{ode:para:full} with parameters in the next table and velocity
\[
   v = \dfrac{\Delta s}{\Delta t}=\dfrac{2\htissuemin (1-r)}{\Delta t},
\]
where $\Delta s$ is the size of the compression of the non-woven fabric. It is set to $\htissuemin (1-r)$, where $r$ is the ratio of the fabric after compression and the presence of the term $2$ is due to the fact that the two fabric sheets are bounded. Finally $\Delta t$ is the time spent during the compression.

\begin{table}[!htb]
\caption{Parameters used in the numerical simulations}
\label{tab:numeric:parameters}
\begin{center}
  \begin{tabular}{|ccl|}
    \hline
    Parameter & Value & Meaning \\
    \hline
    $T_0$         & $20\;\unit{\celsius}$ & 
    Ambient temperature \\
    $\TmaxL$      & $90\;\unit{\celsius}$ &
    Linear model critical temperature \\
    $\TmaxQ$      & $160\;\unit{\celsius}$ &
    Quadratic model critical temperature \\
    $\Cptissue$   & $1800\;\unit{J/(kg\,K)}$ &
    Fabric heat capacity \\
    $\wtissue$    & $0.0126\;\unit{kg/m^2}$ &
    Heat capacity of the fabric \\
    $\htissuemin$ & $14\,\unit{\mu m}$ & 
    Width of the fabric\\
    $r$           & $0.6$ & Compression ratio \\
    $\Delta t$    & $10\,\unit{ms}$, $1\,\unit{ms}$, $0.1\,\unit{ms}$ &
    Bounding time at various tests\\
    $\Ksteel$     & $50\,\unit{W/(m K)}$ & Thermal capacity of steel\\
    \hline
  \end{tabular}
\end{center}
\end{table}
The choice of $\Delta t$ around $1$ $\unit{ms}$ is justified by parameters of the specific bonding process as discussed in the next section.
\par\medskip\noindent
The (numerically) computed heating is plotted in Figure~\ref{fig:heating}.

\begin{figure}[!htb]
  \begin{center}
    \includegraphics[width=5cm]{figure/heating0-1}
    \includegraphics[width=5cm]{figure/heating0-2}
    \includegraphics[width=5cm]{figure/heating0-3}
  \end{center}
  \caption{
  Heating with and without heat flux for a strain of $0.4$ and a compression time of $10 \ \unit{ms}$, $1 \ \unit{ms}$ and $0.1 \ \unit{ms}$, from left to right
  }
  \label{fig:heating}
\end{figure}


\section{The modelling of the bonding process}
\label{bondingprocess}
In this section, we model the bonding process exploited by the company to thermally bond pieces of non-woven fabric, fitting better the geometry of the steel-made rollers and their rotational velocity. 

\newcommand{\VT}{v_{\mathrm{fabric}}}
\newcommand{\Roller}{R}
\newcommand{\cratio}{r}
\newcommand{\VRot}{\omega}
\newcommand{\droller}{d_r}
\newcommand{\Tmodel}{T}
\newcommand{\Tamb}{T_{\textrm{ambient}}}

\begin{figure}[!htb]
  \begin{center}
    \includegraphics[scale=0.3]{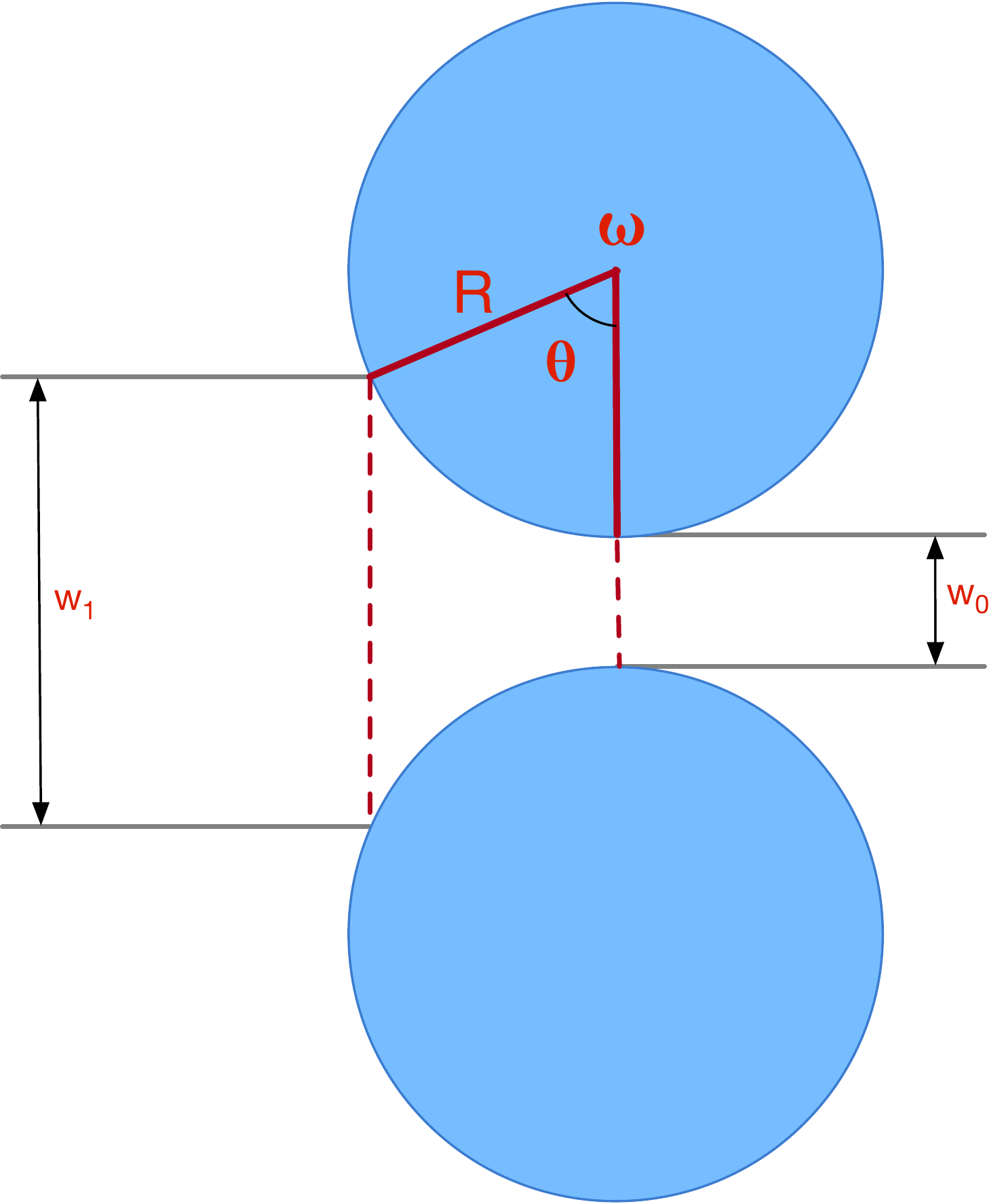}
    \qquad
    \includegraphics[scale=0.3]{figure/tag-geometry.png}
  \end{center}
  \caption{Compression model of two steel-made rollers (left); shape and dimensions of the bonding tags (right) along the roller's lateral surface}
  \label{fig:compression}
\end{figure}

The aim of fast compression is to reach a bonding temperature close to $150 \ \unit{\celsius}$
as noticed by other research~\cite{Hegde2008,Bhat2004} for polypropylene.

\subsection{Determination of the rollers velocity}
Given the velocity of the non-woven fabric on the assembly line $\VT$ (about $360\,\unit{m/min}$, that is $6\,\unit{m/s}$) and the ray of the rollers $\Roller$ (about $0.2 \ \unit{m}$), we can obtain the angular velocity of the rolls:
\begin{equation}\label{eq:omega:def}
    \VRot = \dfrac{\VT}{\Roller} = 
    \dfrac{6 \unit{m/s}}{0.2 \ \unit{m}} = 30
    \unit{\dfrac{rad}{s}}.
\end{equation}

Figure~\ref{fig:compression} (right) shows that a tag, during the roller's rotation, 
spends a time
\[
\Delta t = \dfrac{3.17\,\unit{mm}}{6\,\unit{m/s}} \approx 0.528\,\unit{ms}
\]
at the point of maximum compression. Therefore, the duration of the bonding process on a single tag is on the order of $1\,\unit{ms}$.

\subsection{Determination of the compression angle}
Looking at Figure~\ref{fig:compression} by convention, we assume that the compression starts as the non-woven fabric is compacted (that is $\htissuemin=14\,\unit{\mu m}$) neglecting the effects of the first compaction. Using $\droller$, the distance between the rollers, we define by
\[
  \cratio = \droller/\htissuemin
\]
The compression ratio of the fabric as it passes through the rollers. We can compute the angle $\theta_0$ that the roller has to be subject to in order to take the non-woven fabric from thickness $\htissuemin$ to $\droller$ by solving the following:
\[
    \htissuemin = \droller + 2\Roller(1-\cos\theta_0) = \htissuemin \cratio + 2\Roller(1-\cos\theta_0),
\]
from which
\[
   \cos\theta_0 = 1- \dfrac{\htissuemin(1-\cratio)}{2\Roller}\qquad \textrm{or}
   \qquad
   \theta_0 = \arccos\left( 1- \dfrac{\htissuemin(1-\cratio)}{2\Roller}\right).
\]
Assuming that the angle is much small, we can Taylor approximate $\cos\theta\approx 1-\theta^2/2$:
\begin{equation}\label{eq:theta0}
   \color{blue}
   \theta_0 \approx \dfrac{\sqrt{\htissuemin}}{\sqrt{\Roller}}\sqrt{1-\cratio}.
\end{equation}
For instance, using a ratio $\cratio=0.7$ (that is the non-woven fabric is compressed at $70\%$ w.r.t. its starting size with $\Roller=0.2\,\unit{m}$ and $\htissuemin=14\,\unit{\mu m}$), we obtain
\[
   \theta_0 \approx \dfrac{\sqrt{14\cdot10^{-6}}}{\sqrt{0.2}}\sqrt{1-0.7}\approx 0.004583 \,\unit{rad} = 0.2625\,\unit{deg}.
\]

\subsection{Determination of the bonding time and the velocity profile}
As the compression starts with $\theta(0) = -\theta_0$, the angle is changing with the law
\begin{equation}\label{eq:theta:t}
   \theta(t) = \VRot t - \theta_0
\end{equation}
and the thickness of the non-woven fabric does the same as follows (by Taylor approximating $\cos \theta(t)$ for small angles)
\[
    w(t) = \htissuemin \cratio + 2\Roller(1-\cos \theta(t))\quad\Rightarrow\ \textrm{[Taylor]}\ \Rightarrow\quad
    w(t) := \htissuemin \cratio + \Roller\,\theta(t)^2,
\]
with compression velocity $-w'(t)$
\[
  -w'(t) = -2\Roller\theta(t)\theta'(t)
         = 2\Roller\left(\theta_0-\VRot t\right)\VRot > 0.
\]
The time to perform the bonding is given by $\VRot \Delta t = \theta_0$, which for $\theta_0=0.004583 \,\unit{rad}$ and $\VRot=30 \ \unit{rad/s}$
gives the solution
\[
    \Delta t  = \dfrac{\theta_0}{\VRot} =  \dfrac{0.004583}{30} = 0.00015275\,\unit{s} =  0.15275\,\unit{ms}.
\]

This is the time an infinitesimal section of the fabric spends under compression during the compression phase, based on the parameters used. This duration may vary depending on whether the velocity of the fabric on the roller is increased or decreased. 

Moreover, this time approximately doubles if the release phase is also considered, during which the section remains under compression while the pressure is being released.

\subsection{A more accurate model for heating the non-woven fabric}
\begin{figure}[!htb]
  \begin{center}
    \hspace{-1cm}
    \includegraphics[width=9cm]{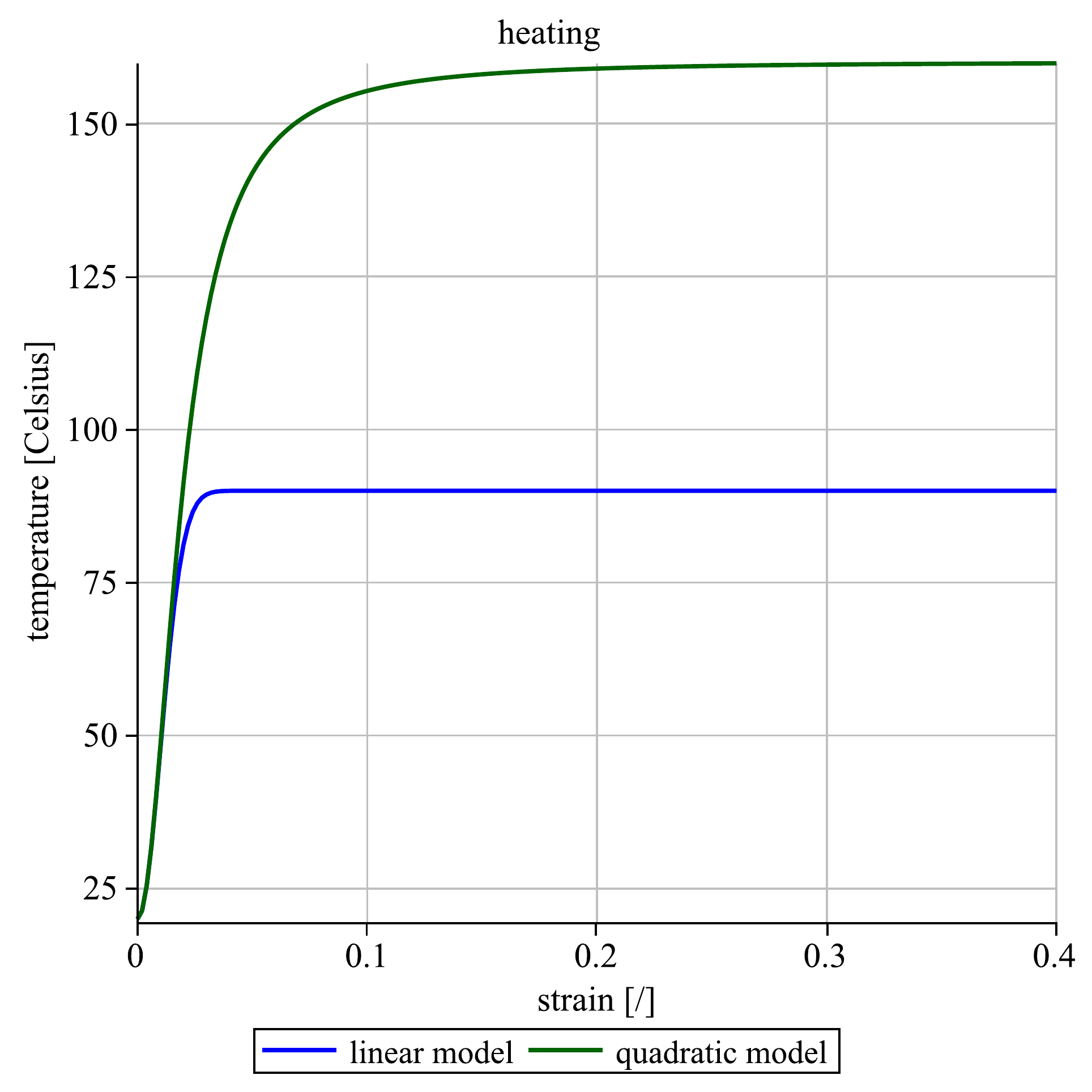}\hspace{-0.5cm}
    \includegraphics[width=9cm]{figure/heating2}
    \hspace{-1cm}\\
    \hspace{-1cm}
    \includegraphics[width=9cm]{figure/heating3}\hspace{-0.5cm}
    \includegraphics[width=9cm]{figure/heating4}
    \hspace{-1cm}
  \end{center}
  \caption{
  Heating for a strain of $0.4$ and different 
  compression time and different thermal conductivity}
  \label{fig:heating2}
\end{figure}

By $w(t)$ and Taylor approximating the $\cos$ ($\cos\theta\approx 1-\theta^2/2$), we can obtain the strain of the fabric w.r.t. the time
$s(t)$
\begin{equation}\label{eq:s:by:time}
   s(t) = 1-\dfrac{w(t)}{\htissuemin} = 
   1-\dfrac{\htissuemin \cratio + \Roller\,\theta(t)^2}{\htissuemin}=
   \color{blue}
   1-\cratio - \dfrac{\Roller}{\htissuemin}(\VRot t - \theta_0)^2,
\end{equation}
from which 
\begin{equation}\label{eq:vel:ratio}
   v(t)=s'(t)=\color{blue}\dfrac{2\Roller\VRot}{\htissuemin}\left(\theta_0-\VRot t\right).
\end{equation}
From \eqref{ode:para}, we obtain the contribution for the heating due to compression.
\begin{equation}\label{eq:model:A}
   \dfrac{\mathrm{d} T(t)}{\mathrm{d}t}  =  \dfrac{\mathrm{d} T(s)}{\mathrm{d}s} s'(t)  =
    \dfrac{\mathrm{d} T(s)}{\mathrm{d}s} v(t)
   =
   \dfrac{v(t)s(t)\Ytissue(s(t))}{\Cptissue\cdot \wtissue}\cdot 
   \max\left(0,\dfrac{\TmaxQ-T(t)}{\TmaxQ-\Tamb}\right)^2,
\end{equation}
where $s(t)$ is given by~\eqref{eq:s:by:time} and $\Ytissue(s)$ by~\eqref{eq:Y:module:tissue}. By~\eqref{eq:T:cool:by:steel}, we have the cooling law due to contact with the steel rolls:
\begin{equation}\label{eq:model:B}
  \dfrac{\mathrm{d}\Tmodel(t)}{\mathrm{d}t}=
  \dfrac{4\Ksteel(\Tsteel-\Tmodel(t))}{\htissuemin(1-s(t))\Cptissue\wtissue}.
\end{equation}
Combining the contributions of~\eqref{eq:model:A} and~\eqref{eq:model:B}, we have the final model
\begin{equation}\label{eq:model:t}
   \dfrac{\mathrm{d} T(t)}{\mathrm{d}t}=
   \dfrac{
    v(t)s(t)\Ytissue(s(t))\cdot 
   \max\left(0,\dfrac{\TmaxQ-T(t)}{\TmaxQ-\Tamb}\right)^2
   +
   \dfrac{4\Ksteel(\Tsteel-\Tmodel(t))}{\htissuemin(1-s(t))}
   }{\Cptissue\wtissue}.
\end{equation}
In order to better compare the solutions, it is suitable to rewrite the equations in terms of the scaled time $\tau$
\[
    t(\tau) = \tau \Delta t = \tau\dfrac{\theta_0}{\VRot}
\]
in such a way that
\begin{equation}\label{eq:time:scale:derivative}
   \dfrac{\mathrm{d} T(\tau)}{\mathrm{d}\tau}=
   \dfrac{\mathrm{d} T(t(\tau))}{\mathrm{d}t}\dfrac{\mathrm{d} t(\tau)}{\mathrm{d}\tau}=
   \dfrac{\mathrm{d} T(t)}{\mathrm{d}t}\dfrac{\theta_0}{\VRot}.
\end{equation}
Moreover, by~\eqref{eq:theta0} and~\eqref{eq:omega:def},
\begin{equation}\label{eq:time:scale}
  \Delta t = 
  \dfrac{\theta_0}{\VRot} 
  =\dfrac{\sqrt{\Roller(1-r)\htissuemin}}{\VT}.
\end{equation}
By~\eqref{eq:theta:t}, \eqref{eq:s:by:time} and~\eqref{eq:vel:ratio} together with~\eqref{eq:theta0}, which gives $R/\htissuemin=(1-r)/\theta_0^2$, we have
\begin{equation}\label{eq:st:tau}
  \begin{split}
   \theta(\tau) = & \VRot t - \theta_0 = (\tau-1)\theta_0,
   \\
   s(\tau)
   = &
   1-\cratio - \dfrac{\Roller}{\htissuemin}(\tau-1)^2\theta_0^2
   = 
   (1-\cratio)(1 - (\tau-1)^2)
   = 
   (1-\cratio)\tau(2-\tau),
   \\
   v(\tau) = &
   \dfrac{2\Roller\VRot}{\htissuemin}(1-\tau)\theta_0=
   \dfrac{2}{\theta_0}\VRot(1-r)(1-\tau),
   \end{split}
\end{equation}
from which $s(0)=0$, $s(1)=1-\cratio$ and $1-s(\tau) = 1-(1-r)\tau(2-\tau)$.
Moreover $v(0)=2\sqrt{\Roller(1-\cratio)/\htissuemin}$ and $v(1)=0$.
Using~\eqref{eq:time:scale} with~\eqref{eq:time:scale:derivative} and~\eqref{eq:model:t}, we obtain
$$
   \color{blue}
   \dfrac{\mathrm{d} T(\tau)}{\mathrm{d}\tau}=
   \dfrac{\theta_0}{\VRot}
   \dfrac{
    v\,s\,\Ytissue(s)\cdot 
   \max\left(0,\dfrac{\TmaxQ-T(\tau)}{\TmaxQ-\Tamb}\right)^2
   +
   \dfrac{4\Ksteel(\Tsteel-\Tmodel(\tau))}{\htissuemin(1-s)}
   }{\Cptissue\wtissue}.
$$
Figure~\ref{fig:heating2} show some possible solutions.

\subsection{Parabolic model}
In order to enhance the compression of the fusing process, we consider a section (along $z$) of the two pieces of non-woven fabric which pass through the rollers (seen as a unique piece of non-woven fabric which is compressed). Neglecting the heat diffusion along the directions $x$ and $y$ (the horizontal ones), since we are assuming that the temperature gradient along such directions is small, we obtain the parabolic equation of the temperature variation (heat equation)
$$
 \wtissue\Cptissue \dfrac{\partial T(t,z)}{\partial t}
 = 
 \Ktissue
 \dfrac{\partial^2 T(t,z)}{\partial^2 z}+\mathcal{Q} (t,z),
$$
for $z\in (-h(t),h(t))$, where $h(t)=(1-s(t))\htissuemin$, and $s(t)$ is the strain as a function of the time given by \eqref{eq:vel:ratio}; the boundary conditions are given by
$$
  \begin{cases}
  T(0,z) = \Tamb, & z\in[-\htissuemin,\htissuemin] \\[1em]
  T(t,-h(t)) = T(t,h(t)) = \Tsteel, & t\in[0,\Delta t] \\[1em]
  \dfrac{\partial T(t,-h(t))}{\partial x} =
  \dfrac{\partial T(t,h(t))}{\partial x}= 0, & t\geq \Delta t
  \end{cases}
$$
where $\Delta t$ is given by \eqref{eq:time:scale}.
The production of the heat $\mathcal{Q} (t,z)$ is due to the pressure on the non-woven fabric and, as a function of the strain, is given by the r.h.s. of \eqref{eq:model:A} scaled on the thickness of the fabric $2h(t)=2(1-s(s))\htissuemin$ and supposed as homogeneous (is not depending on $z$) along the thickness:
$$
  \mathcal{Q}(t,z)=
  \begin{cases}
  \dfrac{v(t)s(t)\Ytissue(s(t))}{2h(t)}\max\left(0,\dfrac{\TmaxQ-T(t,z)}{\TmaxQ-\Tamb}\right)^2,
  & t \leq \Delta t
  \\
  0, & t > \Delta t
  \end{cases}
$$
where, by \eqref{eq:s:by:time} and \eqref{eq:vel:ratio}, we have
$$
  s(t) =
  \begin{cases}
  1-\cratio - \dfrac{\Roller}{\htissuemin}(\VRot t - \theta_0)^2, & t\in[0,\Delta t] \\[1em]
  1-\cratio, &  t\geq \Delta t
  \end{cases} 
  \qquad
  v(t) =
  \dfrac{\Roller}{\htissuemin}
  \begin{cases}
  \theta_0-\VRot t, & t\in[0,\Delta t] \\[1em]
  0, &  t\geq \Delta t.
  \end{cases}
$$
Changing the time in a scaled instant $t = \tau\Delta t$ and $z\in[-h(t),h(t)]$ in $[-1,1]$, i.e. $z=\zeta h(t)$, we obtain the re-scaled equation
$$
 h(\tau)\wtissue\Cptissue \dfrac{\partial T(\tau,\zeta)}{\partial\tau}
 = 
 \Delta t \left(\Ktissue\dfrac{\partial^2 T(\tau,\zeta)}{\partial^2 \zeta}+\overline{\mathcal{Q}}(\tau,\zeta)\right),
$$
where (using \eqref{eq:st:tau}), we have
\[
  \begin{split}
  \overline{\mathcal{Q}}(\tau,\zeta)=&
  \begin{cases}
  \dfrac{v(\tau)s(\tau)\Ytissue(s(\tau))}{2}\max\left(0,\dfrac{\TmaxQ-T(\tau,\zeta)}{\TmaxQ-\Tamb}\right)^2,
  & \tau \leq 1 \\
  0, & \tau > 1
  \end{cases}
  \\
   s(\tau)
   = & (1-\cratio)\tau(2-\tau), \\
   v(\tau) = & \dfrac{2}{\theta_0}\VRot(1-r)(1-\tau), \\
   h(\tau) = &
   \htissuemin \left(1-(1-r)\tau(2-\tau)\right),
  \end{split}
\]
and the boundary conditions become
$$
  \begin{cases}
  T(0,\zeta) = \Tamb, & z\in[-1,1] \\[1em]
  T(\tau,-1) = T(\tau,1) = \Tsteel, & \tau\in[0,1] \\[1em]
  \dfrac{\partial T(\tau,-1)}{\partial \zeta} =
  \dfrac{\partial T(\tau,1)}{\partial \zeta}= 0, & \tau\geq 1.
  \end{cases}
$$
Using a spatial discretization  $\Delta\zeta=2/N$ and $\zeta_k=-1+k\Delta\zeta$, we sample the temperature along the points $\zeta_k$ with the functions $T_k(\tau)\approx T(\tau,\zeta_k)$ for which, after the discretization, we obtain an ODE system for $\tau \leq 1$
$$
  \left\{
  \begin{aligned}
    T_0'(\tau)
    =\,& 0,
    \\
    T_k'(\tau)
    =\,& \Delta t\dfrac{
    \big(\Ktissue/(h(\tau)\Delta\zeta^2)\big) 
    \big(T_{k-1}(\tau)-2T_k(\tau)+T_{k+1}(\tau)\big)
    +\overline{\mathcal{Q}}_k(\tau)
    }{\wtissue\Cptissue},
    \\
    T_N'(\tau)
    =\,& 0,
    \\
    \overline{\mathcal{Q}}_k(\tau)=&
    \dfrac{v(\tau)s(\tau)\Ytissue(s(\tau))}{2}
    \max\left(0,\dfrac{\TmaxQ-T_k(\tau)}{\TmaxQ-\Tamb}\right)^2,
  \end{aligned}
  \right.
$$
and one for $\tau\geq 0$
$$
  \left\{
  \begin{aligned}
    T_0'(\tau)
    =& 
    \dfrac{\Delta t\Ktissue}{\Delta\zeta^2h(\tau)\wtissue\Cptissue}
    \big(T_1(\tau)-T_0(\tau)\big),
    \\
    T_k'(\tau)
    =& 
    \dfrac{\Delta t\Ktissue}{\Delta\zeta^2h(\tau)\wtissue\Cptissue}
    \big(T_{k-1}(\tau)-2T_k(\tau)+T_{k+1}(\tau)\big),
    \\
    T_N'(\tau)
    =&
    \dfrac{\Delta t\Ktissue}{\Delta\zeta^2h(\tau)\wtissue\Cptissue}
    \big(T_{N-1}(\tau)-T_N(\tau)\big).
  \end{aligned}
  \right.
$$

\begin{figure}[!hb]
  \begin{center}
    \includegraphics[scale=0.3]{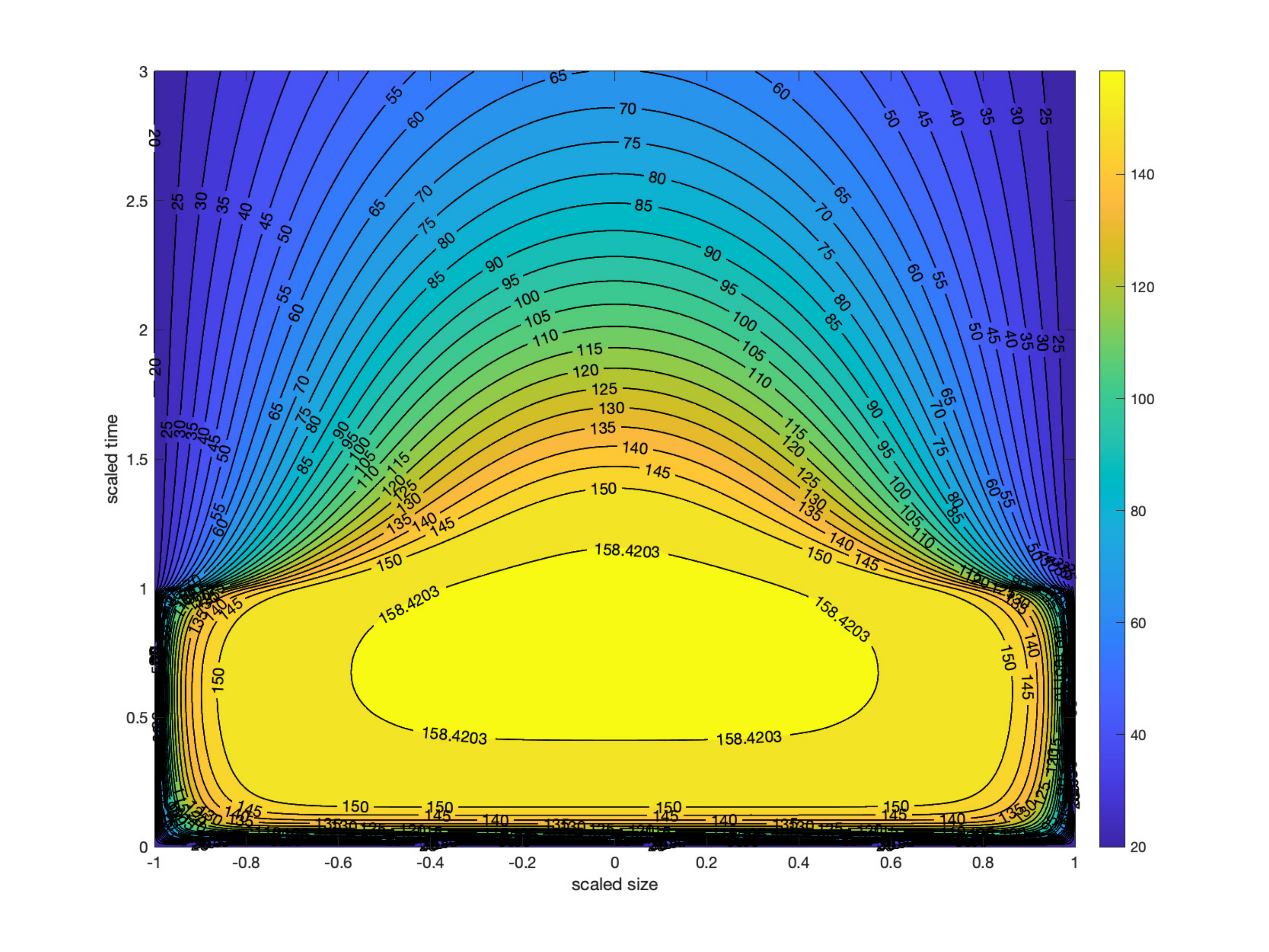}
    \includegraphics[scale=0.3]{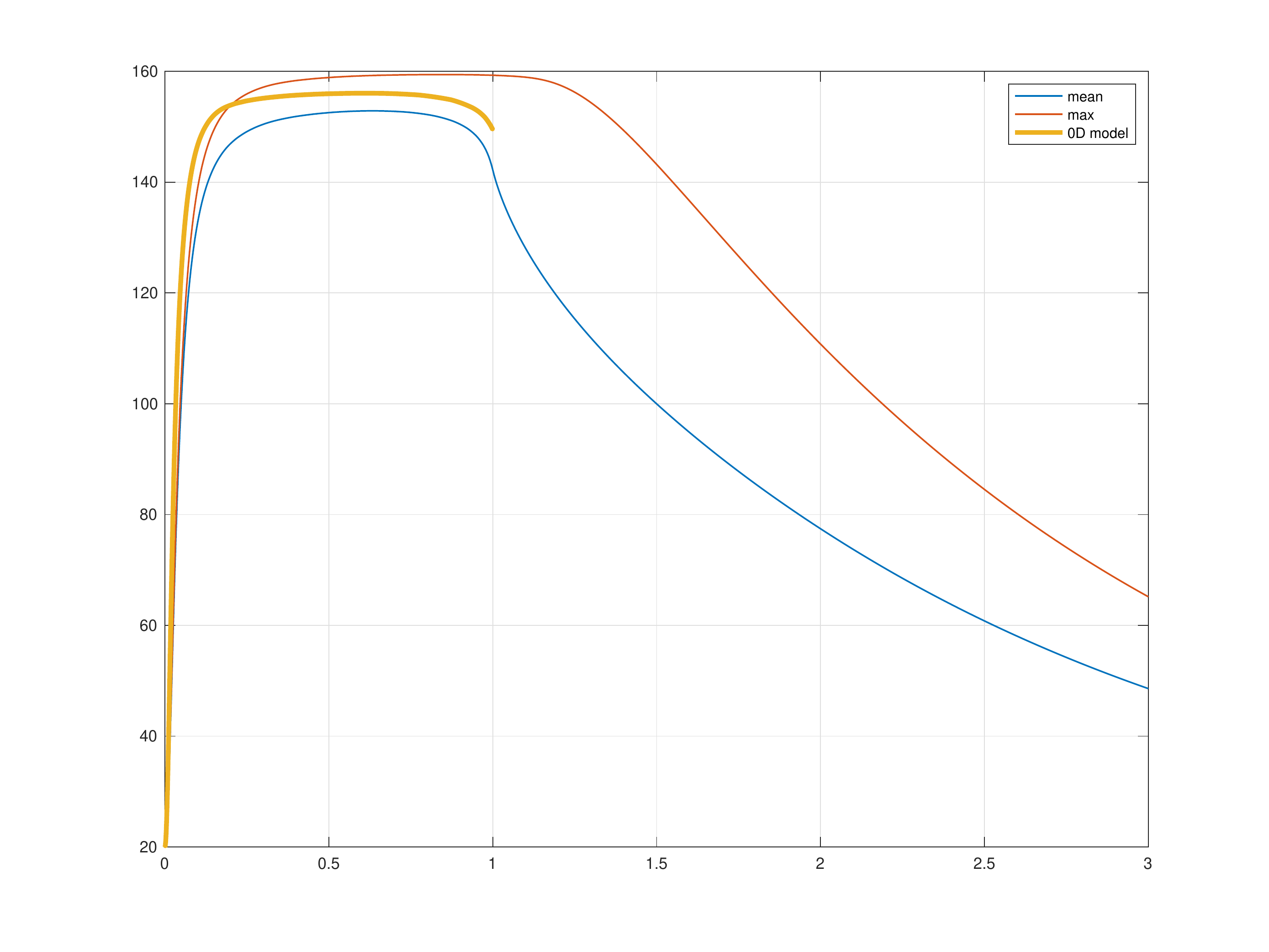}
    \includegraphics[scale=0.6]{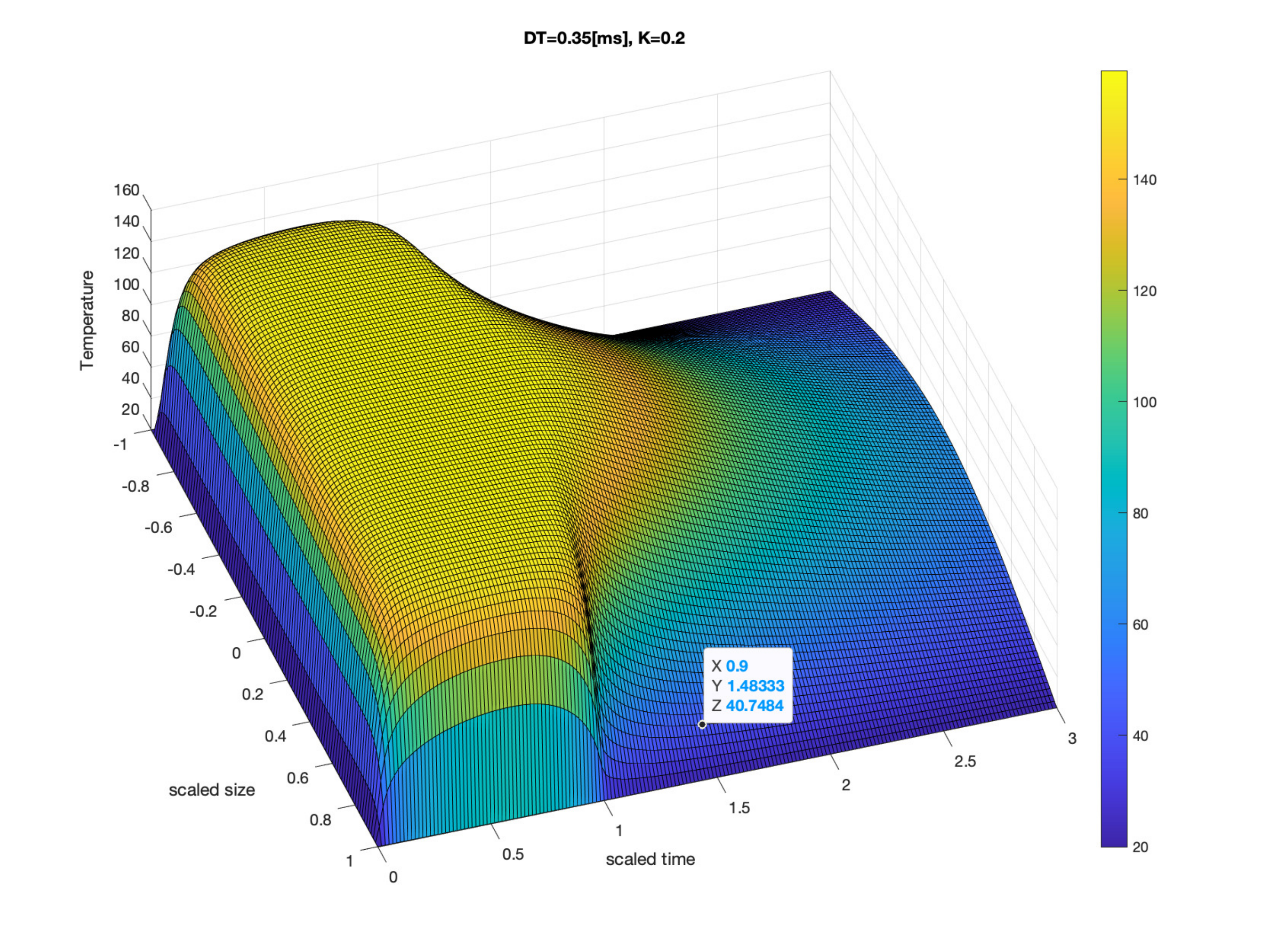}
  \end{center}
  \caption{Solution with the parabolic model and standard parameters:
  $\Ksteel=17\,\unit{W}/(\unit{m K})$, $\VT=6\,\unit{m/s}$ and $r=0.8$.
  Estimated outgoing homogeneous temperature of $25-55 \ \unit{\celsius}$.
 A large yellow area indicates effective bonding.
A wide extent along the scaled coordinate suggests that a significant portion of the fabric is melted.
A broad region with elevated temperature implies that the heat is sustained long enough to ensure proper fusion of the fabric.
  }
  \label{fig:para1}
\end{figure}

\begin{figure}[!hb]
  \begin{center}
    \includegraphics[scale=0.3]{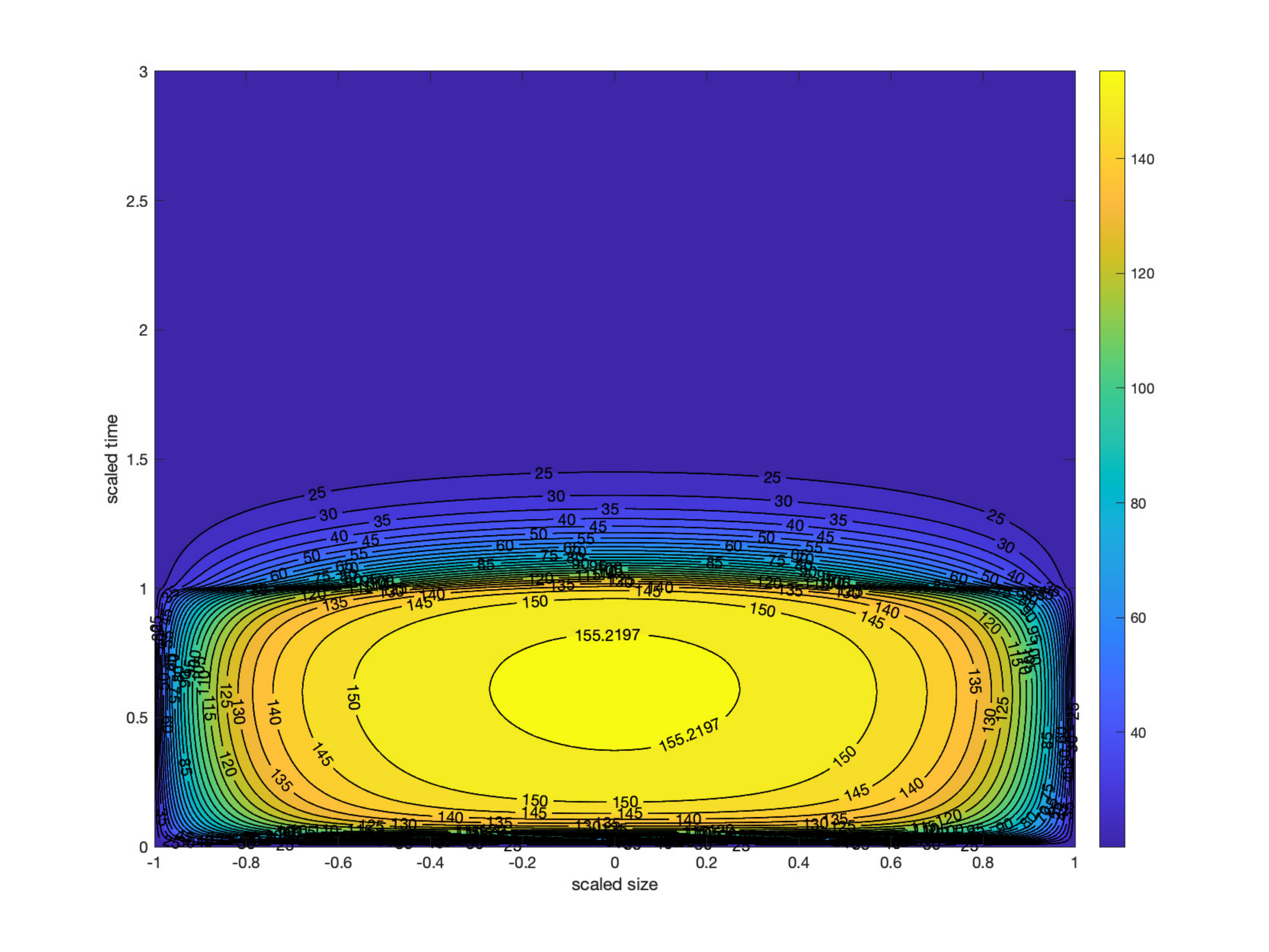}
    \includegraphics[scale=0.3]{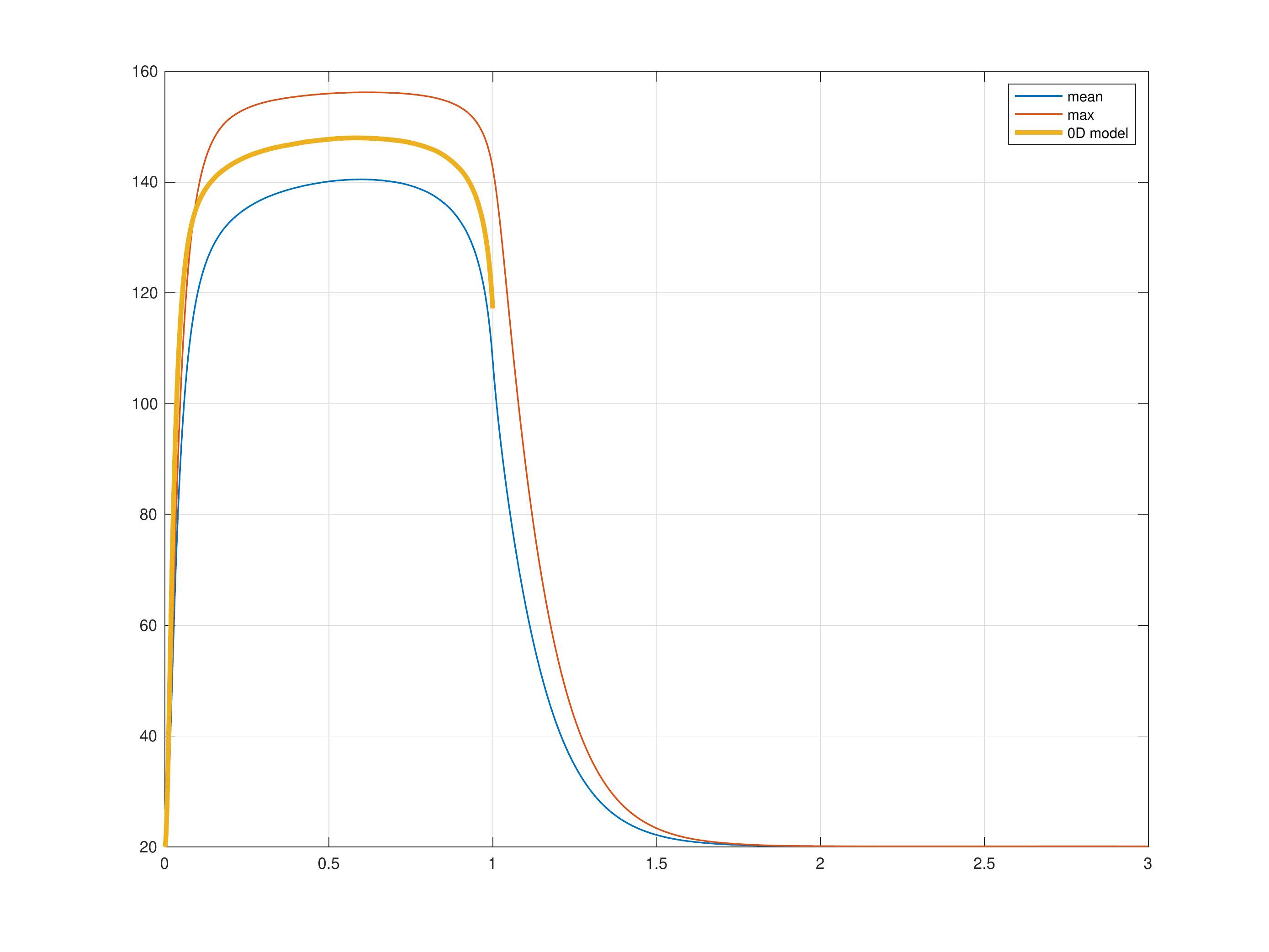}
    \includegraphics[scale=0.6]{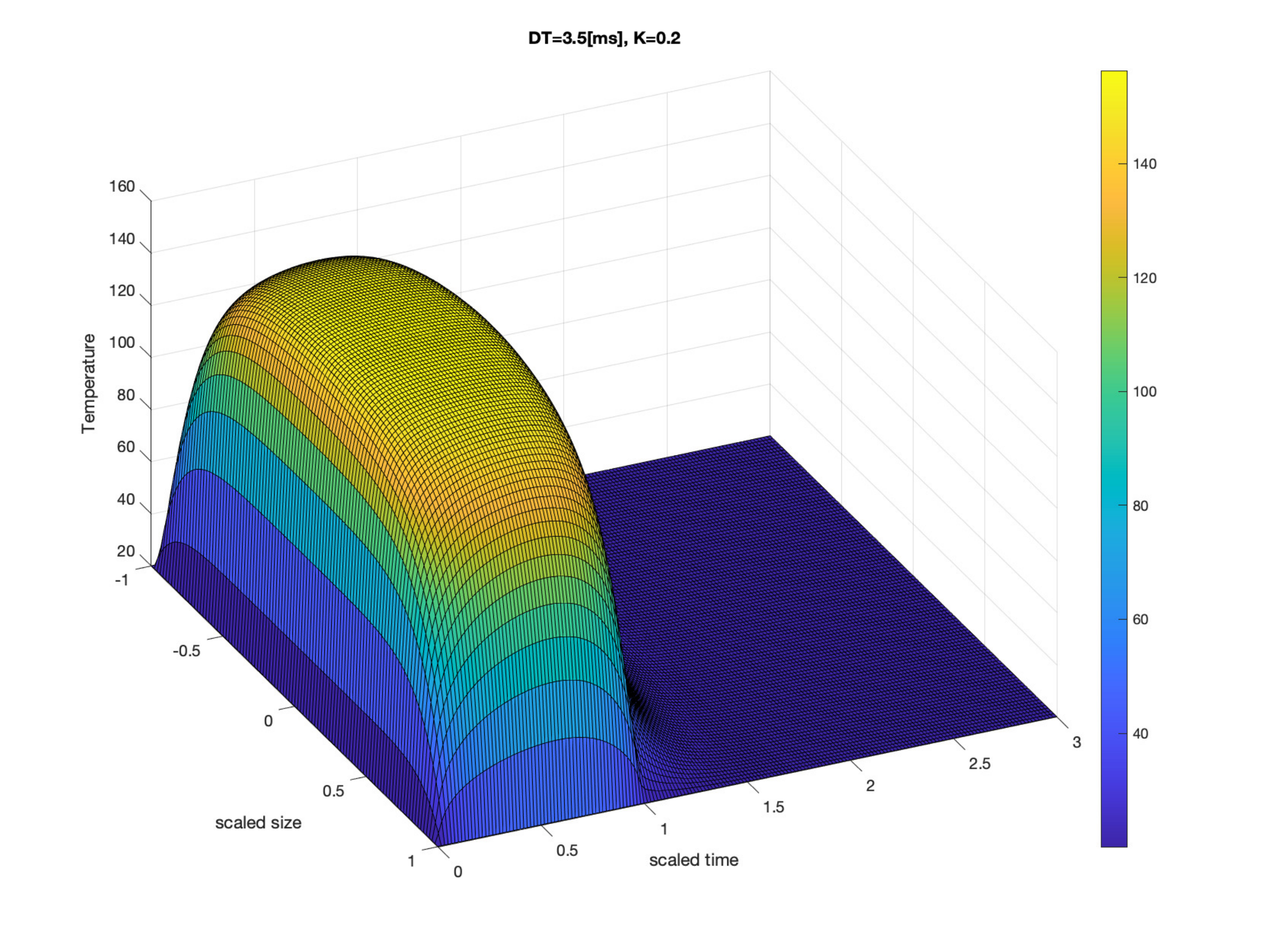}
  \end{center}
  \caption{Solution with the parabolic model and slow scrolling:
  $\Ksteel=17\,\unit{W}/(\unit{m K})$, $\VT=0.6\,\unit{m/s}$ and $r=0.8$.
  The peak temperature is localized within a narrow region, indicating that the resulting bond is likely to be weak or entirely absent.
  A good bound is obtained with flat high temperature in space and time.
  }
  \label{fig:para2}
\end{figure}

\begin{figure}[!hb]
  \begin{center}
    \includegraphics[scale=0.3]{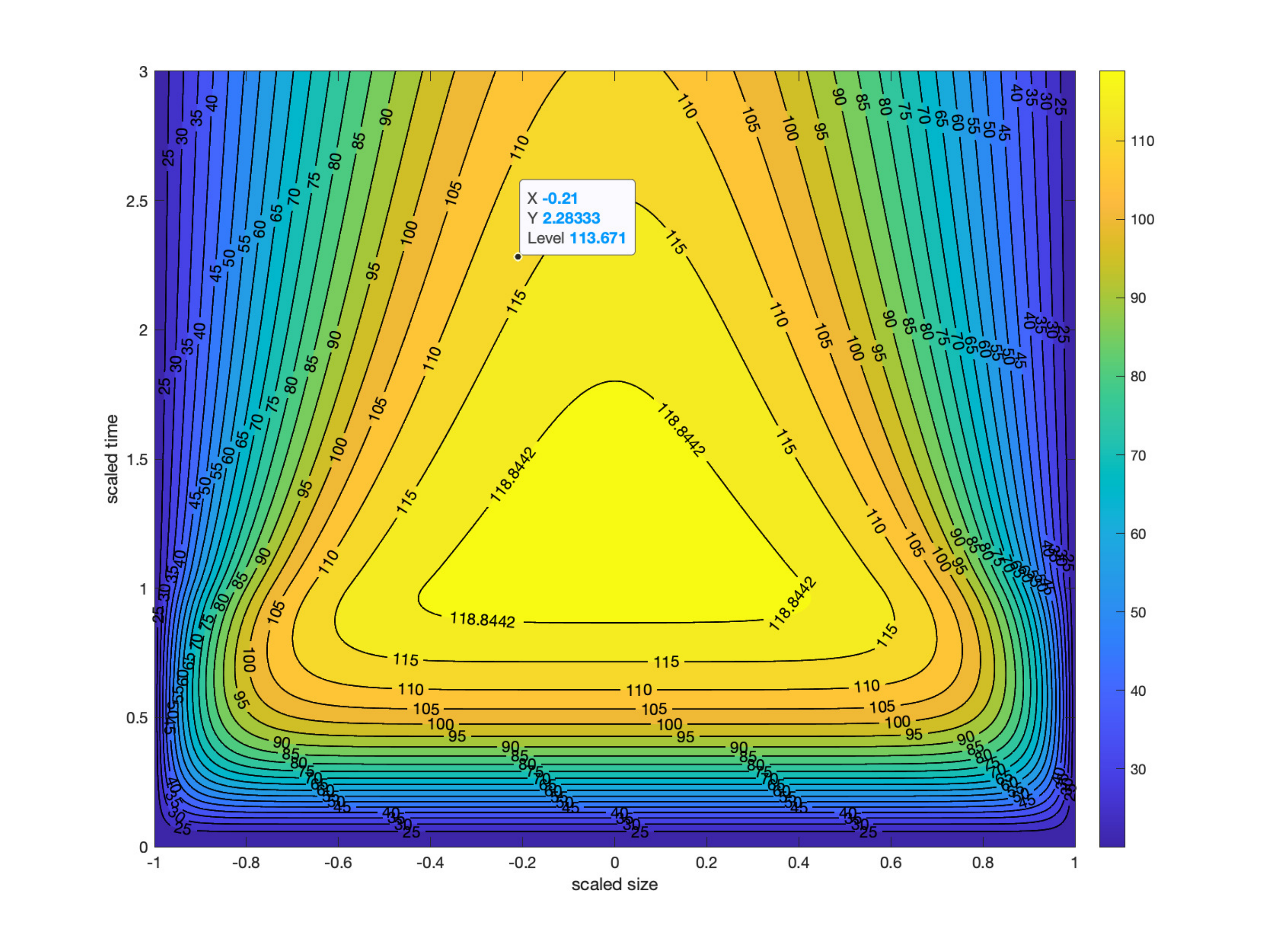}
    \includegraphics[scale=0.3]{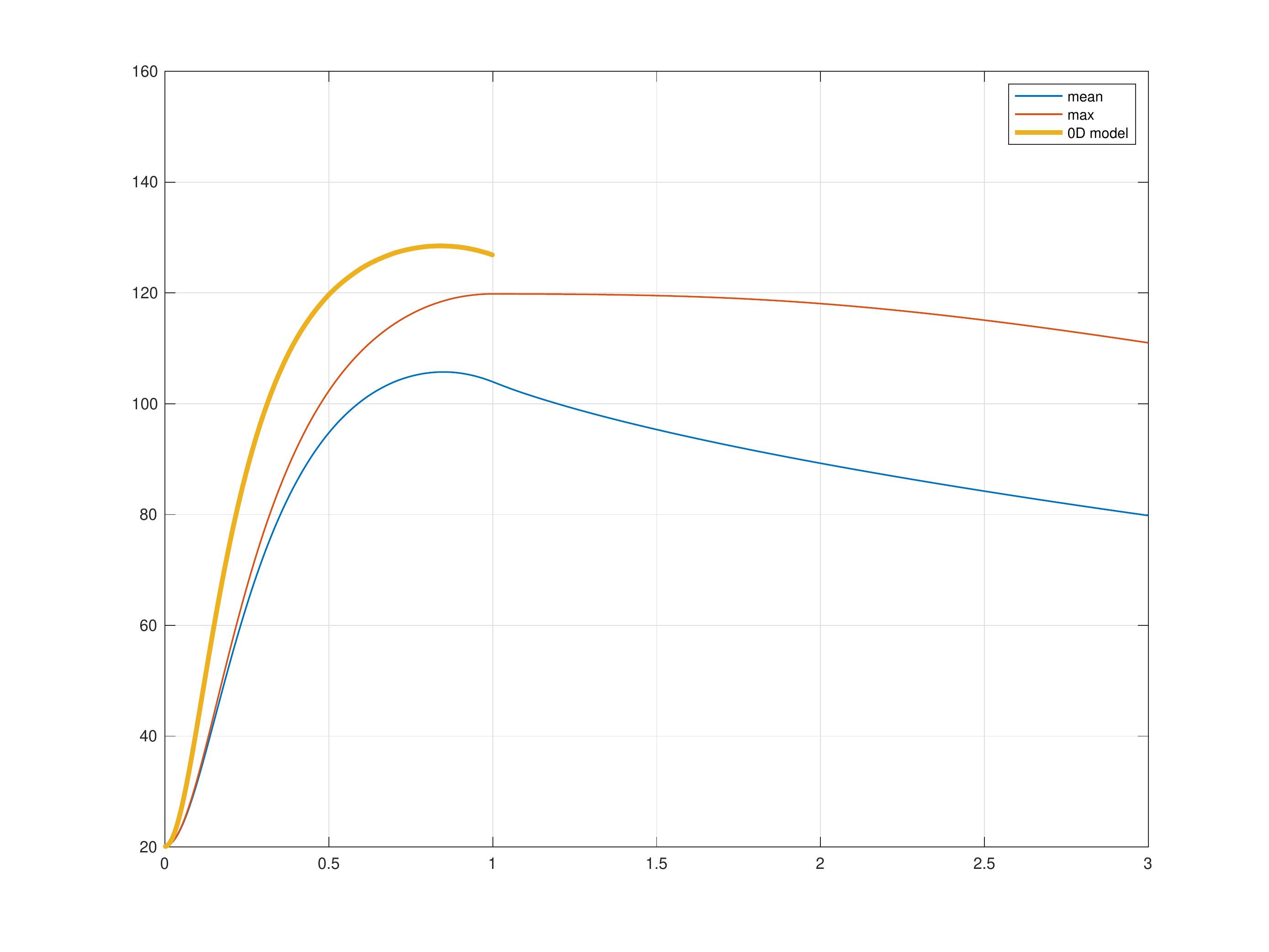}
    \includegraphics[scale=0.6]{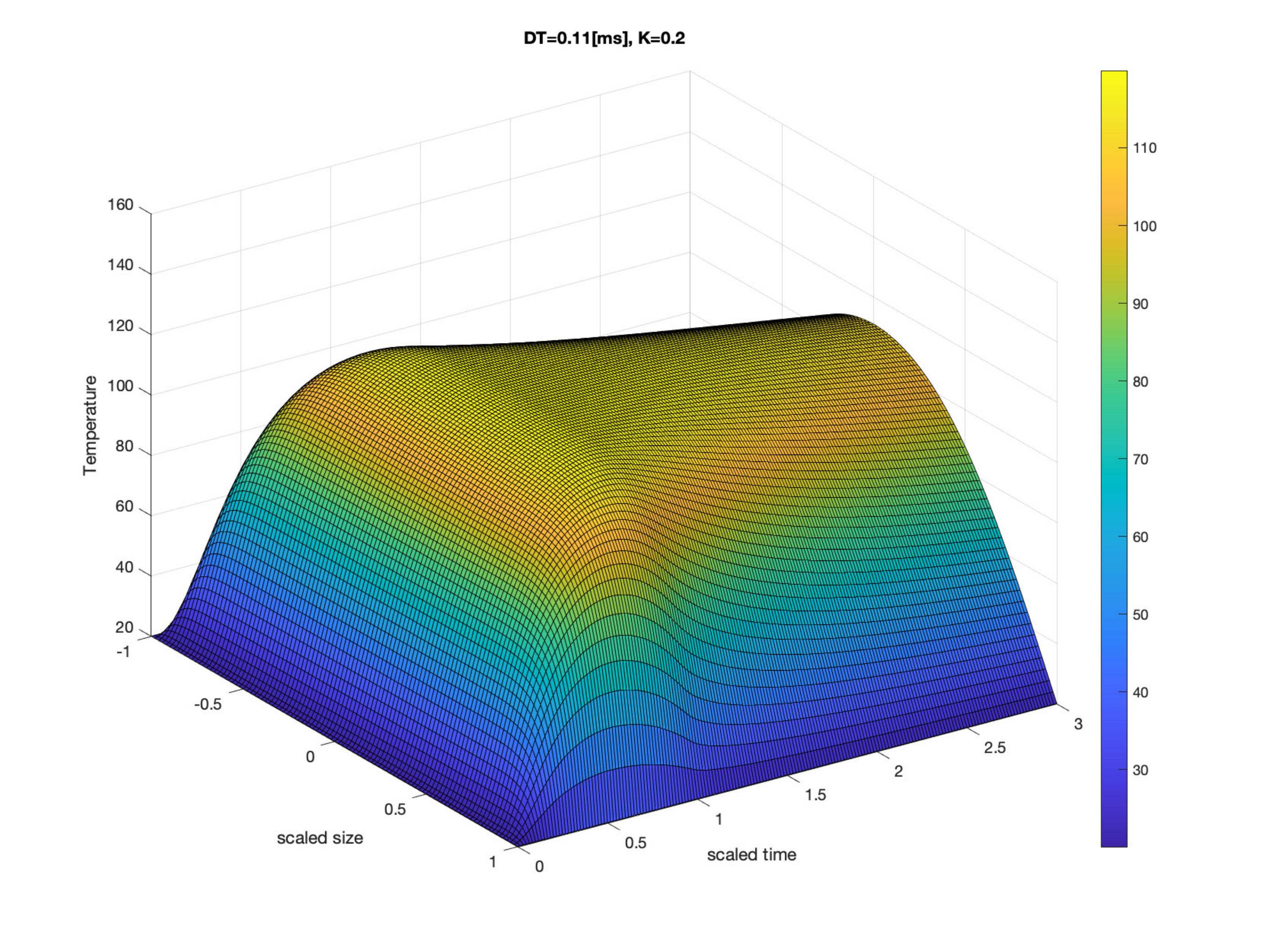}
  \end{center}
  \caption{Solution with the parabolic model with lower compression:
  $\Ksteel=17\,\unit{W}/(\unit{m K})$, $\VT=6\,\unit{m/s}$ and  $r=0.95$.
  The peak temperature is too low, indicating that the resulting bond is likely to be absent.
  }
  \label{fig:para3}
\end{figure}

Figures~\ref{fig:para1}, \ref{fig:para2} 
and \ref{fig:para3} shows the results 
of the numerical simulations with 
varius settings.

\section{Conclusions}
\label{conclusions}

We have proposed and studied a simplified model of the bonding process based on a limited and measurable set of parameters related to the textile fabric and the bonding machinery.
All of these parameters are explicitly incorporated into the analytical framework. Textile fabric parameters can be easily estimated or are typically available from technical data sheets.
The machinery parameters -- such as compression speed, compression ratio, and ambient temperature -- can be readily configured on-site.

Given the parameters of a specific fabric, numerical simulations allow for an evaluation of how different settings in the machinery influence the efficiency of the bonding process. This, in turn, enables optimization of the bonding conditions in various textile fabrics. As a result, expensive machine downtimes typically required for empirical adjustments may be significantly reduced or even replaced by preliminary simulation-based analyses.

For simulations, we used the line method, following a spatial discretization of the governing partial differential equation (PDE) using finite difference techniques~\cite{Formaggia}.

\begin{itemize}
\item \emph{0D} simplified models provide a first-order approximation of the temperature evolution and exposure time, offering a preliminary assessment of the feasibility of achieving effective bonding.

\item \emph{2D}
models enable more accurate predictions of bond quality based on both fabric and process parameters; thus, it is possible to estimate the performance of the final product in different scenarios.
\end{itemize}

The main contributions of this study can be summarized as follows:

\begin{itemize}

	\item A reduced-order model of the bonding process has been developed, relying on a limited number of parameters that are measurable or readily available in data sheets.
	
    \item The model supports the identification of optimal machinery settings, such as compression speed, compression ratio, and ambient temperature, facilitating faster setup and fine-tuning in industrial environments.
	
    \item The model offers predictive insight into bond quality based on input parameters. Furthermore, it can determine whether bonding is theoretically feasible (e.g., whether the required fusion temperature is reached and maintained long enough for bonding to occur).
	
    \item Whether a new fabric is introduced, the model simplifies the reconfiguration of the bonding process. It enables assessment of bonding feasibility under the current setup (e.g., compression rollers), and allows for adjustments -- such as modifying compression profiles -- to extend process applicability and reduce downtime.

\end{itemize}

Future developments will focus on constructing and analyzing a mechanical-analytical model of the fibers behavior before and after the bonding process. This model will incorporate the elastoplastic nature of the fibers, represented through suitable mechanical analogs such as combinations of springs and dampers.








\clearpage

\baselineskip=0.9\normalbaselineskip
\phantomsection\addcontentsline{toc}{section}{\numberline{}References}

\bibliographystyle{CAIM_Sciendo_bibstyle}

\bibliography{main}

\end{document}